\documentclass[12pt]{amsart}
\usepackage{amsmath,amsfonts,amssymb,amsthm}
\usepackage{hyperref}
\newtheorem{Theorem}{Theorem}[section]
\newtheorem{Lemma}[Theorem]{Lemma}
\newtheorem{Corollary}[Theorem]{Corollary}
\newtheorem{Proposition}[Theorem]{Proposition}
\theoremstyle{definition}

\newtheorem{Example}[Theorem]{Example}
\newtheorem{Remark}[Theorem]{Remark}
\newtheorem{Problem}[Theorem]{Problem}

\newcommand{\C}{\mathbb C}
\newcommand{\N}{\mathbb N}
\newcommand{\ve}{\varepsilon}
\newcommand{\ov}{\overline}
\newcommand{\pa}{\partial}
\newcommand{\fr}{\frac}
\newcommand{\Sect}{\mathcal S}
\newcommand{\Lev}{\mathcal L}
\newcommand{\Gam}{\Gamma_{\!f}}
\DeclareMathOperator{\Arg}{Arg}
\DeclareMathOperator{\conv}{conv}
\renewcommand{\Re}{\operatorname{Re}}
\renewcommand{\Im}{\operatorname{Im}}

\title[Uniform algebras generated by powers and smooth functions]{On the uniform algebras generated by powers of the coordinates and smooth functions on compact subsets of $\C^N$}
\author{Kieu Phuong Chi and Le Minh Thien}
\address{Kieu Phuong Chi, \newline
Department of Mathematics and Applications, Saigon University, \newline
273 An Duong Vuong, Ho Chi Minh City, Vietnam}
\email{kieuphuongchi@sgu.edu.vn}
\address{Le Minh Thien
\newline
\indent Department of Mathematics and Applications,
\newline
\indent Saigon  University, 273  An Duong Vuong, Ho Chi Minh City, Vietnam}
\email{lmthien.lqd@hcm.edu.vn}


\begin{document}

\begin{abstract}
Let $K\subset\C^N$ be compact and let $f=(f_1,\dots,f_N)$ be a smooth map on a neighbourhood of $K$. We give sufficient conditions for the algebra generated by $z_1^{m_1},\dots,z_N^{m_N}$ and $f_1,\dots,f_N$ (or by all monomials of degree $n$ and $f_1,\dots,f_N$) to be dense in $C(K)$. We appropriately extend the techniques developed in our previous works \cite{Chi,CT,DC1}  and \cite{DC}  for algebras of functions on the complex plane to tackle the problem in the several-variable setting. More precisely, the argument separates a purely algebraic reduction, in which the rotated graphs are glued one coordinate at a time, from the study of the single graph of $f$. For the latter we give an exact computation of the Levi form of $|w-f|^{2p}$, which yields a Hermitian matrix criterion on the polynomial hull of the graph; a criterion by peak functions that removes the invertibility hypothesis from a theorem of Preskenis; a division theorem for $f_j=\bar k_jh_j$ with degenerate holomorphic factors $h_j$; and a treatment of the singular set $\{\det(\pa f_i/\pa\bar z_j)=0\}$, which in several variables is not finite. The angle $\pi/m_j$ in the reduction is shown to be sharp.
\end{abstract}

\maketitle

\section{Introduction and statement of the main results}

For compact $K\subset\C^N$ and $g_1,\dots,g_k\in C(K)$ let $P[g_1,\dots,g_k;K]$ be the uniform closure on $K$ of the polynomials in $g_1,\dots,g_k$. The study of $P[z,f;D]$ for a disk $D$ goes back to Wermer \cite{W}; Preskenis \cite{Pre} treated orientation reversing $f$ and the higher dimensional analogue on the polydisk; Minsker \cite{M} proved $P[z^m,\bar z^n;D]=C(D)$ for coprime $m,n$; and de Paepe \cite{P3,PW}, Dieu--de Paepe \cite{DP}, Dieu--Chi \cite{DC},Chi \cite{Chi} and Bharali \cite{B} developed the local theory of $P[z^m,f;D]$ on small disks. The paper \cite{CT} considers the global question for $P[z^n,f;K]$ with $K\subset\C$ polynomially convex. Here we take $K\subset\C^N$.

For $m=(m_1,\dots,m_N)\in\N^N$ write $z^{[m]}=(z_1^{m_1},\dots,z_N^{m_N})$. Let $\Gam(A)=\{(z,f(z)):z\in A\}$, $\Gam=\Gam(K)$, and for $\theta>0$ let $\Sect_\theta=\{s\in\C\setminus\{0\}:|\arg s|<\theta\}$. A prototype of our results reads as follows (see Theorems \ref{thmA}, \ref{thmHess}).

\medskip
\noindent\textbf{Theorem.} \emph{Let $K\subset\C^N$ be compact and polynomially convex, $f\in C^2$ near $K$, and $c_j\in\N$ with $\gcd(c_j,m_j)=1$. Assume:
\begin{enumerate}
\item[(i)] $P(\Gam(Z))=C(\Gam(Z))$, where $Z=\{a\in K:\det(\pa f_i/\pa\bar z_j)(a)=0\}$;
\item[(ii)] $z_j^{c_j}f_j(z)\in\Sect_{\pi/m_j}$ whenever $m_j\ge2$, $z\in K$, $z_j\neq0$;
\item[(iii)] for some integer $p\ge1$ the Hermitian form $\Lev_{f,p}(z,w;\cdot)$ of \eqref{eq:form} is positive semidefinite at every point $(z,w)$ of the polynomial hull of $\Gam$.
\end{enumerate}
Then $P[z_1^{m_1},\dots,z_N^{m_N},f_1,\dots,f_N;K]=C(K)$.}
\medskip

Compared with $N=1$ there are three new features.

(1) The set $Z$ is the zero set of one complex function on $\C^N$, so it has real codimension two in general and cannot be assumed finite: for $f=(\bar z_1^2,\bar z_2^2)$ it is $\{z_1z_2=0\}$. Section \ref{sec:base} treats it by stratification.

(2) The preimage of $\{(z^{[m]},f(z))\}$ consists of $\prod m_j$ graphs, and two of them meet over $K\cap\{z_j=0,\ j\in J\}$, which is not a point. Theorem \ref{thmA} avoids any analysis of these intersections by rotating one coordinate at a time; the sector conditions only involve the $j$-th coordinate at the $j$-th stage. If instead one keeps the scalar action of the roots of unity, the pieces meet at one point exactly as for $N=1$, and this leads to the second algebra, generated by the monomials of degree $n$ (Theorem \ref{thmB}).

(3) Duval's theorem \cite{Du} has no available global analogue. Its local estimate carries over (Proposition \ref{prop:locDuval}), and we replace the global statement by Proposition \ref{prop:peak}, Theorem \ref{thmMon} and the division result Proposition \ref{prop:div}. The latter extends to all dimensions a one-variable phenomenon: $f_j=\bar k_jh_j$ with $h_j$ allowed to vanish on $\{k_j=0\}$, without any differential inequality (Theorem \ref{thmTwist}).
\par
In Section \ref{sec:HW} we use the theorem of H\"ormander--Wermer, Weinstock and Browder--Cole--Wermer: $P[z,\bar z+R;K]=C(K)$ for any compact $K$ and $R$  Lipschitz with constant $<1$ on $K$ (\cite[Thm.~1.6.10]{St}, Theorem \ref{thmLip}). Combined with the reduction, this gives density of $P[z^{[m]},f;K]$ and $P[\nu_n,f;K]$ on compact sets that are not polynomially convex (Corollary \ref{corHW}, Example \ref{ex:shell}). We also discuss stability under $C^2$ perturbations (Corollary \ref{cor:stab}, Remark \ref{rem:instab}).
\par In Section \ref{sec:one} we return to the case $N=1$ and, through a series of examples, establish connections with several known results
\section{Preliminaries}

Let $\widehat X$ be the polynomial hull of a compact $X\subset\C^N$, that is
$$\widehat{X}=\{z\in \Bbb C^N: \vert p(z)\vert\leq \max_{x\in X}\vert p(x)\vert,\;\text{for every holomorphic polynomial}\; p\; \text{on}\; \C^N
\}.$$
 Let $P(X)$ be the uniform closure of polynomials on $X$, and $\ov\Delta$ the closed unit disk. Recall: the projection of $\widehat X$ to a factor lies in the hull of the projection of $X$; if $X$ is polynomially convex, it has a neighbourhood basis of Runge domains $D$ (interiors of polynomial polyhedra), every function holomorphic near $X$ belongs to $P(X)$ (Oka--Weil), and $P(X)=C(X)$ implies polynomial convexity. If $x\in\widehat X$ there is a probability measure $\mu$ on $X$ with $g(x)=\int g\,d\mu$ for all $g\in P(X)$, and $\mu$ can be chosen Jensen (\cite{Ga}), i.e.\ $\log|g(x)|\le\int\log|g|\,d\mu$.

If $g\colon K\to\C^L$ is continuous and injective on $K$, then $P[g_1,\dots,g_L;K]=C(K)$ if and only if $P(g(K))=C(g(K))$. In particular $P[z,f;K]=C(K)$ iff $P(\Gam)=C(\Gam)$.

If $K$ is polynomially convex then $\widehat{\Gam}\subset K\times\C^N$; so if moreover $\Gam$ is polynomially convex and $A\subset K$ is polynomially convex, then $\Gam(A)$ is polynomially convex, because $\widehat{\Gam(A)}\subset\Gam\cap(A\times\C^N)$. \hfill(2.1)

\begin{Theorem}[{\cite[Thm.~4.3.4]{Ho}, \cite[Thm.~1.3.11]{St}}]\label{thm:psh}
Let $D\subset\C^M$ be a Runge domain, $X\subset D$ compact with $\widehat X\subset D$, and $\varphi$ plurisubharmonic on $D$. Then $\varphi\le\sup_X\varphi$ on $\widehat X$.
\end{Theorem}

The next lemma contains all we need of Kallin's lemma and of its version by Stout and de Paepe.

\begin{Lemma}\label{lem:zero}
Let $X=X_1\cup X_2$ be compact with $P(X_i)=C(X_i)$, $i=1,2$. If $g\in P(X)$ vanishes on $X_1$ and has no zero on $X_2\setminus X_1$, then $P(X)=C(X)$.
\end{Lemma}
\begin{proof}
By the Hahn--Banach and Riesz theorems, it suffices to show that a complex regular Borel measure $\mu$ on $X$ with $\int p\,d\mu=0$ for all polynomials $p$ vanishes. Such a $\mu$ annihilates $P(X)$. Since $P(X)$ is an algebra containing $g$, the measure $\nu=g\mu$ also annihilates $P(X)$, and in particular all polynomials. As $g=0$ on $X_1$, $\nu$ is carried by $X_2\setminus X_1\subset X_2$. Since $P(X_2)=C(X_2)$, a measure on $X_2$ annihilating the polynomials vanishes, so $\nu=0$. As $g\ne0$ on $X_2\setminus X_1$, we get $|\mu|(X_2\setminus X_1)=0$. Thus $\mu$ is a measure on $X_1$ annihilating the polynomials, and $\mu=0$ since $P(X_1)=C(X_1)$.
\end{proof}

\begin{Corollary}[{Kallin--Stout--de Paepe \cite{Ka,P3,St}}]\label{cor:stout}
Let $X_1,X_2\subset\C^M$ be compact with $P(X_i)=C(X_i)$. Let $Y_1,Y_2\subset\C$ be compact with $Y_1\cap Y_2=\{0\}$, $Y=Y_1\cup Y_2$ having connected complement and $0\notin\operatorname{int}Y$. If $p$ is a polynomial with $p(X_i)\subset Y_i$ and $p^{-1}(0)\cap(X_1\cup X_2)=X_1\cap X_2$, then $P(X_1\cup X_2)=C(X_1\cup X_2)$.
\end{Corollary}
\begin{proof}
Let $h$ be the function on $Y$ equal to $0$ on $Y_1$ and to $s$ on $Y_2$. It is well defined, as $Y_1\cap Y_2=\{0\}$ and $h(0)=0$ in both definitions, and it is continuous since $Y_1,Y_2$ are closed. As $0\notin\operatorname{int}Y$, every component of $\operatorname{int}Y$ is a connected subset of $Y\setminus\{0\}=(Y_1\setminus\{0\})\sqcup(Y_2\setminus\{0\})$, a disjoint union of two relatively open sets. So it lies in one of them, and $h$ is holomorphic on it. By Mergelyan's theorem, $h$ is a uniform limit on $Y$ of polynomials $r_l$. Then $g=h\circ p$ is the uniform limit on $X=X_1\cup X_2$ of the polynomials $r_l\circ p$, so $g\in P(X)$. It vanishes on $X_1$. On $X_2\setminus X_1$ it equals $p$, which has no zero there because $p^{-1}(0)\cap X=X_1\cap X_2$. Lemma \ref{lem:zero} applies.
\end{proof}

Let $G\subset GL_N(\C)$ be a finite group. It acts on $\C^N\times\C^M$ by $g(\zeta,w)=(g\zeta,w)$. For $f\colon K\to\C^M$ continuous, let $L_g=\{(gz,f(z)):z\in K\}$ and $L=\bigcup_{g\in G}L_g$. The following replaces the lemma of de Paepe (\cite{P3}, \cite[Lemma 2.5]{Chi}).

\begin{Lemma}\label{lem:avg}
Let $\mathcal I$ be a finite set of generators of the invariant ring $\C[z]^G$. If $\mathcal I\cup\{f_1,\dots,f_M\}$ separates the points of $K$ and $P(L)=C(L)$, then $P[\mathcal I,f_1,\dots,f_M;K]=C(K)$.
\end{Lemma}
\begin{proof}
Let $\phi\in C(K)$ and define $\tilde\phi$ on $L$ by $\tilde\phi(gz,f(z))=\phi(z)$.

\emph{Well defined.} Let $(gz,f(z))=(g'z',f(z'))$. Then $z'=g'^{-1}gz$ lies in the $G$-orbit of $z$, so every element of $\mathcal I$ takes the same value at $z$ and $z'$. Also $f_i(z)=f_i(z')$ for all $i$. By the separation hypothesis $z=z'$.

\emph{Continuous and invariant.} Each $L_g$ is compact, and $z\mapsto(gz,f(z))$ is a continuous injection of the compact set $K$ onto it (injective because $g$ is), hence a homeomorphism. So $\tilde\phi$ is continuous on each of the finitely many closed sets $L_g$, hence on $L$. Further $\tilde\phi(g'\zeta,w)=\tilde\phi(\zeta,w)$ on $L$ for $g'\in G$, by the definition of $\tilde\phi$.

\emph{Averaging.} Take a polynomial $q$ on $\C^N\times\C^M$ with $|q-\tilde\phi|<\ve$ on $L$, and put $q^G(\zeta,w)=|G|^{-1}\sum_{g}q(g\zeta,w)$. For $(\zeta,w)\in L$ we have $(g\zeta,w)\in L$ and $\tilde\phi(g\zeta,w)=\tilde\phi(\zeta,w)$, so
$$|q^G(\zeta,w)-\tilde\phi(\zeta,w)|\le|G|^{-1}\sum_g|q(g\zeta,w)-\tilde\phi(g\zeta,w)|<\ve.$$
Expand $q^G=\sum_\beta c_\beta(\zeta)w^\beta$. Each $c_\beta$ is $G$-invariant, so it is a polynomial in $\mathcal I$. Thus $q^G(\zeta,w)=Q(\mathcal I(\zeta),w)$ for a polynomial $Q$. Restricting to $L_1=\Gam$ gives $|Q(\mathcal I(z),f(z))-\phi(z)|<\ve$ on $K$.
\end{proof}

A real $C^1$ submanifold is totally real at a point if its tangent space there contains no complex line. For the approximation step we use the following form of a theorem of O'Farrell--Preskenis--Walsh, Stout and Izzo.

\begin{Theorem}[{\cite{OPW}, \cite[Thm.~6.3.2, Cor.~6.3.3]{St}, \cite[Thm.~2.1]{Iz}}]\label{thm:Iz}
Let $X\subset\C^M$ be compact and polynomially convex and let $E\subset X$ be compact (possibly empty) with $P(E)=C(E)$. If every point of $X\setminus E$ has a neighbourhood in which $X$ is contained in a $C^1$ totally real submanifold, then $P(X)=C(X)$.
\end{Theorem}
In \cite{St} the hypothesis is that $X\setminus E$ is a totally real subset of $\C^M\setminus E$. The proof there is local: it shows that a measure orthogonal to the functions holomorphic near $X$ vanishes near each point of $X\setminus E$, using only that $X$ lies locally in a $C^1$ totally real manifold. In all our applications $X\setminus E$ lies in one $C^1$ totally real manifold, namely a graph over a complex manifold, which is the form of \cite[Cor.~6.3.3]{St}.

For $f\in C^1$ put $A_f=(\pa f_i/\pa z_j)$ and $B_f=(\pa f_i/\pa\bar z_j)$.

\begin{Lemma}\label{lem:tr}
Let $S$ be a complex submanifold of an open subset of $\C^N$, $f\colon S\to\C^M$ the restriction of a $C^1$ map, and $a\in S$. Then $\{(z,f(z)):z\in S\}$ is totally real at $(a,f(a))$ if and only if $B_f(a)\bar v\neq0$ for every $v\in T_aS\setminus\{0\}$.
\end{Lemma}
\begin{proof}
Write $A=A_f(a)$, $B=B_f(a)$. The real differential of $f$ at $a$ is $v\mapsto Av+B\bar v$, so the tangent space of the graph at $(a,f(a))$ is $\mathcal T=\{\tau(v):=(v,Av+B\bar v):\ v\in T_aS\}$. It contains a complex line if and only if there is $v\ne0$ with $i\tau(v)\in\mathcal T$. Now $i\tau(v)=(iv,iAv+iB\bar v)$. Comparing first components, $i\tau(v)\in\mathcal T$ if and only if $i\tau(v)=\tau(iv)$. Here $iv\in T_aS$ as $S$ is complex, and
$$\tau(iv)=\big(iv,\ A(iv)+B\ov{iv}\big)=(iv,\ iAv-iB\bar v).$$
So $i\tau(v)=\tau(iv)$ if and only if $2iB\bar v=0$, that is, $B\bar v=0$.
\end{proof}
For $S$ open in $\C^N$, $M=N$, this is Chirka's and Preskenis' condition $\det B_f(a)\ne0$ (\cite{Ch}, \cite[Thm.~3.2]{Pre}).

\begin{Lemma}\label{lem:change}
Let $K\subset\C^N$ be polynomially convex and $\Omega\supset K$ a Runge domain. Let $\Lambda\colon\Omega\to GL_M(\C)$ and $b\colon\Omega\to\C^M$ be holomorphic, and $\Phi(z,w)=(z,\Lambda(z)w+b(z))$ on $D=\Omega\times\C^M$. For compact $X\subset K\times\C^M$: $X$ is polynomially convex iff $\Phi(X)$ is; in that case $P(X)=C(X)$ iff $P(\Phi(X))=C(\Phi(X))$.
\end{Lemma}
\begin{proof}
The map $\Phi$ is a biholomorphism of $D$, with inverse $\Phi^{-1}(z,w)=(z,\Lambda(z)^{-1}(w-b(z)))$, and it maps $K\times\C^M$ onto itself. Further $D$ is a Runge domain: expanding in powers of $w$ and approximating the coefficients on $\Omega$, every function in $\mathcal O(D)$ is a limit of polynomials, uniformly on compact subsets of $D$.

For a compact $Y\subset K\times\C^M$ we have $\widehat Y\subset K\times\C^M\subset D$, because the projection of $\widehat Y$ to $\C^N$ lies in $\widehat K=K$. Let $\widehat Y_D$ be the hull of $Y$ with respect to $\mathcal O(D)$. Then $\widehat Y_D\subset\widehat Y$, as polynomials belong to $\mathcal O(D)$. Conversely, if $y\in\widehat Y$ and $g\in\mathcal O(D)$, then $g$ is a uniform limit of polynomials on the compact subset $\widehat Y$ of $D$, so $|g(y)|\le\sup_Y|g|$. Hence $\widehat Y=\widehat Y_D$.

Since $g\mapsto g\circ\Phi$ is a bijection of $\mathcal O(D)$, we have $\Phi(\widehat X_D)=\widehat{\Phi(X)}_D$. With the above, $\Phi(\widehat X)=\widehat{\Phi(X)}$, which proves the first claim. Let $X$ be polynomially convex. The closure in $C(X)$ of $\mathcal O(D)|_X$ is $P(X)$, again as $D$ is Runge, and likewise for $\Phi(X)$. The isometric isomorphism $C(\Phi(X))\to C(X)$, $\phi\mapsto\phi\circ\Phi$, carries $\mathcal O(D)|_{\Phi(X)}$ onto $\mathcal O(D)|_X$, hence $P(\Phi(X))$ onto $P(X)$. This proves the second claim.
\end{proof}

\section{Reduction to a single graph}

\begin{Lemma}[Rotation gluing]\label{lem:rot}
Let $E\subset\C^M$ be compact with $P(E)=C(E)$, $T\in GL_M(\C)$ with $T^n=I$, $n\ge2$, and $q$ a polynomial with $q\circ T=\omega q$, where $\omega$ is a primitive $n$-th root of unity. Assume $q(E)\subset\{0\}\cup\Sect_{\pi/n}$ and that $Tx=x$ for every $x\in E$ with $q(x)=0$. Then $P(\bigcup_{s=0}^{n-1}T^sE)=C(\bigcup_{s=0}^{n-1}T^sE)$.
\end{Lemma}
\begin{proof}
Put $E_s=T^sE$ and $Y_s=\omega^s\conv(q(E)\cup\{0\})$. Since $n\ge2$, the set $\Sect_{\pi/n}\cup\{0\}$ is convex, so $Y_s\subset\omega^s(\Sect_{\pi/n}\cup\{0\})$; as $\omega^s$ runs through all $n$-th roots of unity, the $Y_s$ are compact, convex and meet pairwise exactly at $0$, and $q(E_s)=\omega^sq(E)\subset Y_s$. Every union of the $Y_s$ is compact and star-shaped about $0$, hence has connected complement, and does not contain a neighbourhood of $0$, since the rays $\arg s=\pi/n+2\pi t/n$ miss it. If $x=T^sx_0\in E_s$ and $q(x)=0$ then $q(x_0)=0$, so $Tx_0=x_0$ and $x=x_0$. Hence $Q:=E\cap\{q=0\}=E_s\cap\{q=0\}$ for all $s$. Let $X_1=E_0\cup\dots\cup E_{s-1}$, $X_2=E_s$. If $x\in X_1\cap X_2$ then $q(x)\in Y_s\cap\bigcup_{t<s}Y_t=\{0\}$, so $x\in Q$; thus $X_1\cap X_2=Q=q^{-1}(0)\cap(X_1\cup X_2)$. Also $P(E_s)=C(E_s)$, $T$ being linear and invertible. Corollary \ref{cor:stout} and induction on $s$ give the lemma.
\end{proof}

\begin{Theorem}\label{thmA}
Let $K\subset\C^N$ be compact, $f\colon K\to\C^M$ continuous, $m\in\N^N$. For each $j$ with $m_j\ge2$ let $c_j\in\N$, $\gcd(c_j,m_j)=1$, and $h_j$ a polynomial on $\C^N\times\C^M$, and put
$$\Phi_j(z)=z_j^{c_j}\,h_j\big(z^{[m]},f(z)\big).$$
Assume that
\begin{equation}\label{eq:sector}
\Phi_j(z)\in\Sect_{\pi/m_j}\quad\text{for all }z\in K\text{ with }z_j\neq0.
\end{equation}
Then $z_1^{m_1},\dots,z_N^{m_N},f_1,\dots,f_M$ separate the points of $K$, and
$$P\big[z_1^{m_1},\dots,z_N^{m_N},f_1,\dots,f_M;K\big]=C(K)\iff P\big[z_1,\dots,z_N,f_1,\dots,f_M;K\big]=C(K).$$
\end{Theorem}
The default choice is $M=N$, $h_j(s,w)=w_j$, for which \eqref{eq:sector} is condition (ii) of the Introduction.
\begin{proof}
\emph{Separation.} Let $z,z'\in K$ have the same image. Then $z'_j=\rho_jz_j$ with $\rho_j^{m_j}=1$, and $\Phi_j(z')=\rho_j^{c_j}\Phi_j(z)$. If $z_j\ne0$ and $m_j\ge2$, both values lie in the sector of opening $2\pi/m_j$, and they differ by the $m_j$-th root of unity $\rho_j^{c_j}$; hence $\rho_j^{c_j}=1$ and $\rho_j=1$, as $\gcd(c_j,m_j)=1$. So $z=z'$.

\emph{Density.} ``$\Rightarrow$'' is trivial. For ``$\Leftarrow$'', let $G=\prod_j\mu_{m_j}$ act diagonally; $\C[z]^G$ is generated by the $z_j^{m_j}$, and by Lemma \ref{lem:avg} it is enough to prove $P(L)=C(L)$. Reindex so that $m_j\ge2$ exactly for $j\le r$. Let $T_j$ multiply $\zeta_j$ by $e^{2\pi i/m_j}$ and fix the remaining coordinates of $(\zeta,w)$, and put $E_0=\Gam$, $E_j=\bigcup_sT_j^sE_{j-1}$; so $E_r=L$ and $P(E_0)=C(E_0)$ by hypothesis. Let $q_j(\zeta,w)=\zeta_j^{c_j}h_j(\zeta^{[m]},w)$; then $q_j\circ T_j=\omega_jq_j$ with $\omega_j=e^{2\pi ic_j/m_j}$ primitive. A point of $E_{j-1}$ is $x=(gz,f(z))$ with $z\in K$ and $g=\operatorname{diag}(\rho_1,\dots,\rho_{j-1},1,\dots,1)$, so $q_j(x)=\Phi_j(z)$, independent of $g$. Thus $q_j(E_{j-1})\subset\{0\}\cup\Sect_{\pi/m_j}$, and $q_j(x)=0$ forces $z_j=0$, whence $\zeta_j=0$ and $T_jx=x$. Lemma \ref{lem:rot} gives $P(E_j)=C(E_j)$ inductively.
\end{proof}

\begin{Remark}\label{rem:separation}
Separation is a consequence of \eqref{eq:sector}. This also applies to $N=1$: the alternative ``(a) $\Leftrightarrow$ (b)'' in \cite{CT} is really the statement that (b) holds, since (a) follows from the angular condition.
\end{Remark}

Let $\nu_n(z)$ be the vector of all monomials of degree $n$ in $z_1,\dots,z_N$. These generate the invariants of the scalar group $\mu_n\cdot I$.

\begin{Theorem}\label{thmB}
Let $K\subset\C^N$ be compact, $f\colon K\to\C^M$ continuous, $n\ge2$, $\omega=e^{2\pi i/n}$, and let $q$ be a polynomial with $q(\omega\zeta,w)=\omega^cq(\zeta,w)$ for some $c$ coprime to $n$. If $\Phi(z):=q(z,f(z))\in\Sect_{\pi/n}$ for all $z\in K\setminus\{0\}$, then $\nu_n,f$ separate the points of $K$ and
$$P[\nu_n,f_1,\dots,f_M;K]=C(K)\iff P[z_1,\dots,z_N,f_1,\dots,f_M;K]=C(K).$$
The standard choice, for $M=N$, is $q=\sum_j\zeta_jw_j$, i.e.\ $\Phi(z)=\sum_jz_jf_j(z)$.
\end{Theorem}
\begin{proof}
\emph{Separation.} Let $\nu_n,f$ agree at $z,z'\in K$. If $z=0$ then $\nu_n(z')=0$, so $z_j'^{\,n}=0$ for all $j$ and $z'=0$. Let $z\ne0$, so $z'\ne0$. The equality $\nu_n(z')=\nu_n(z)$ says $z'^\alpha=z^\alpha$ for $|\alpha|=n$. Choose $j$ with $z_j\ne0$ and put $\lambda=z_j'/z_j$, so $\lambda^n=1$. Comparing the monomials $z_iz_j^{n-1}$ gives $z_i'=\lambda z_i$ for all $i$. Thus $z'=\omega^tz$ for some $t$. As $f(z')=f(z)$,
$$\Phi(z')=q(\omega^tz,f(z))=\omega^{ct}\Phi(z).$$
Both values lie in $\Sect_{\pi/n}$, an open sector of opening $2\pi/n$, and they differ by an $n$-th root of unity. Hence $\omega^{ct}=1$. As $\gcd(c,n)=1$, $\omega^t=1$ and $z'=z$.

\emph{Density.} ``$\Rightarrow$'' is trivial. Let $P(\Gam)=C(\Gam)$. Apply Lemma \ref{lem:rot} with $E=\Gam$, $T(\zeta,w)=(\omega\zeta,w)$, the polynomial $q$ and the primitive root $\omega^c$. We check the hypotheses. Evaluating $q(\omega\zeta,w)=\omega^cq(\zeta,w)$ at $\zeta=0$ gives $(1-\omega^c)q(0,w)=0$, so $q(0,w)=0$. Thus $q(E)\subset\{0\}\cup\Sect_{\pi/n}$. If $x=(z,f(z))\in E$ and $q(x)=0$, then $z=0$, hence $Tx=x$. So $P(L)=C(L)$ with $L=\bigcup_sT^s\Gam$. The invariants of $\mu_nI$ are generated by the monomials of degree $n$, and by the first part Lemma \ref{lem:avg} applies. Here the pieces $L_k$ meet only over the origin, exactly as in one variable.
\end{proof}

The scalar action can be replaced by any cyclic diagonal action, which gives algebras generated by other sets of monomials.

\begin{Theorem}\label{thmC}
Let $K\subset\C^N$ be compact, $f\colon K\to\C^M$ continuous, $n\ge2$, $\omega=e^{2\pi i/n}$, $a\in\mathbb Z^N$ with $\gcd(a_1,\dots,a_N,n)=1$, $T_a=\operatorname{diag}(\omega^{a_1},\dots,\omega^{a_N})$, and let $\mathcal I_a$ be a finite set of monomials generating $\C[z]^{\langle T_a\rangle}$, that is, the monomials $z^\alpha$ with $\sum_ja_j\alpha_j\equiv0\pmod n$. Let $q$ be a polynomial with $q(T_a\zeta,w)=\omega^cq(\zeta,w)$, $\gcd(c,n)=1$, and put $\Phi(z)=q(z,f(z))$. Assume that $\Phi(K)\subset\{0\}\cup\Sect_{\pi/n}$, and that $\Phi(z)=0$ implies $T_az=z$. Then $\mathcal I_a,f$ separate the points of $K$ and
$$P[\mathcal I_a,f_1,\dots,f_M;K]=C(K)\iff P[z_1,\dots,z_N,f_1,\dots,f_M;K]=C(K).$$
\end{Theorem}
\begin{proof}
Let $\mathcal I_a,f$ agree at $z,z'$. The invariants separate orbits, so $z'=T_a^tz$. If $\Phi(z)=0$ then $z'=z$. Otherwise $\Phi(z')=\omega^{ct}\Phi(z)$ with both values in $\Sect_{\pi/n}$, so $\omega^{ct}=1$, hence $n\mid t$ and $z'=z$. Now apply Lemma \ref{lem:rot} with $E=\Gam$ and $T(\zeta,w)=(T_a\zeta,w)$, and then Lemma \ref{lem:avg} with $G=\langle T_a\rangle$.
\end{proof}
Theorem \ref{thmB} is the case $a=(1,\dots,1)$, in which the condition on the zeros of $\Phi$ is automatic, since $\Phi\ne0$ on $K\setminus\{0\}$.

\begin{Remark}
The same proof applies to a finite group $G\subset GL_N(\C)$ with a chain $G=G_0\rhd G_1\rhd\dots\rhd G_r=1$ with cyclic quotients, given polynomials $q_i$ transforming under $G_{i-1}$ by a character $\chi_i$ with kernel $G_i$, provided that a zero $z$ of $q_i(z,f(z))$ has a stabiliser mapping onto $G_{i-1}/G_i$. For this one uses that Lemma \ref{lem:rot} remains true, with the same proof, if $T^nE=E$ and $Tx\in E$ for every $x\in E$ with $q(x)=0$. Theorems \ref{thmA} and \ref{thmB} are the two main instances.
\end{Remark}

\begin{Remark}[Sharpness of the angle]\label{rem:sharp}
The constant $\pi/m_j$ in \eqref{eq:sector} cannot be increased. Let $m\ge2$, $\theta>\pi/m$, $a=e^{-i\pi/m}$, $\rho=e^{2\pi i/m}$, $K=[0,a]\cup[0,\rho a]\subset\C$ and $f\equiv1$. Then $zf=z\in\Sect_\theta$ on $K\setminus\{0\}$, as $\arg a=-\pi/m$ and $\arg(\rho a)=\pi/m$. By Lavrentiev's theorem $P[z,f;K]=C(K)$. But $a^m=(\rho a)^m$ and $f(a)=f(\rho a)$, so $z^m,f$ do not separate $K$.
\end{Remark}

\section{The single graph}\label{sec:base}

Fix $f\in C^2$ near $K\subset\C^N$ with values in $\C^N$. For $\xi\in\C^N$ write $L_{f_i}(z;\xi)=\sum_{j,k}\fr{\pa^2f_i}{\pa z_j\pa\bar z_k}(z)\xi_j\bar\xi_k$, and $\langle a,b\rangle=\sum a_i\bar b_i$. For $w\ne f(z)$ set $u=w-f(z)$, $e=u/|u|$, and for an integer $p\ge1$
\begin{equation}\label{eq:form}
\Lev_{f,p}(z,w;\xi)=|B_f\bar\xi|^2+\Big(1-\fr1p\Big)\big|\langle B_f\bar\xi,e\rangle\big|^2-2\Re\sum_{i=1}^N\bar u_i\,L_{f_i}(z;\xi).
\end{equation}
For $p=1$ this is the Hermitian form with matrix
\begin{equation}\label{eq:matrix}
\begin{gathered}
M_f(z,w)_{jk}=\sum_{i=1}^N\Big(\ov{\fr{\pa f_i}{\pa\bar z_j}}\,\fr{\pa f_i}{\pa\bar z_k}-\bar u_i\,\fr{\pa^2f_i}{\pa z_j\pa\bar z_k}-u_i\,\ov{\fr{\pa^2f_i}{\pa z_k\pa\bar z_j}}\Big),\\
\Lev_{f,1}(z,w;\xi)=\sum_{j,k}M_f(z,w)_{jk}\,\xi_j\bar\xi_k.
\end{gathered}
\end{equation}
Note the dependence on $w$. When $N=1$,
$$\Lev_{f,p}(z,w;\xi)=\Big[\big(2-\tfrac1p\big)|f_{\bar z}|^2-2\Re\big((\bar w-\bar f)f_{z\bar z}\big)\Big]|\xi|^2.$$

\begin{Proposition}\label{prop:levi}
Let $G_p(z,w)=|w-f(z)|^{2p}$. At a point with $w\ne f(z)$, the function $G_p$ has positive semidefinite Levi form if and only if $\Lev_{f,p}(z,w;\cdot)\ge0$ on $\C^N$. At points with $w=f(z)$ the Levi form of $G_p$ is positive semidefinite.
\end{Proposition}
\begin{proof}
Let $F=|u|^2=\sum_i|h_i|^2$ with $h_i=w_i-f_i$, and take a tangent vector $v=(\xi,\eta)\in\C^N\times\C^N$. Put $\zeta=\eta-A_f\xi$, $B=B_f$. For a $C^2$ complex function $h$ write $\pa h(v)=\sum_jh_jv_j$, $\bar\pa h(\bar v)=\sum_kh_{\bar k}\bar v_k$ and $L_h(v)=\sum_{j,k}h_{j\bar k}v_j\bar v_k$. Differentiating $h\bar h$,
$$\pa_j\pa_{\bar k}(h\bar h)=h_{j\bar k}\bar h+h_{\bar k}\,\ov{h_{\bar j}}+h_j\,\ov{h_k}+h\,\ov{h_{k\bar j}} ,$$
and multiplying by $v_j\bar v_k$ and summing, we get
$$L_{|h|^2}(v)=|\pa h(v)|^2+\big|\bar\pa h(\bar v)\big|^2+2\Re\big(\bar h\,L_h(v)\big).$$
For $h=h_i$ we have $\pa h_i(v)=\eta_i-(A_f\xi)_i=\zeta_i$, $\bar\pa h_i(\bar v)=-(B\bar\xi)_i$ and $L_{h_i}(v)=-L_{f_i}(z;\xi)$. Also $\pa\bar h_i(v)=\ov{\bar\pa h_i(\bar v)}=-\ov{(B\bar\xi)_i}$, so $\pa F(v)=\sum_i(\bar h_i\zeta_i-h_i\ov{(B\bar\xi)_i})$. Summing over $i$ gives
\begin{equation}\label{eq:LF}
L_F(v)=|\zeta|^2+|B\bar\xi|^2-2\Re\sum_i\bar u_iL_{f_i}(z;\xi),\qquad \pa F(v)=\langle\zeta,u\rangle-\langle u,B\bar\xi\rangle.
\end{equation}
Since $G_p=F^p$, the chain rule gives
\begin{equation}\label{eq:chain}
L_{G_p}(v)=pF^{p-1}L_F(v)+p(p-1)F^{p-2}|\pa F(v)|^2 .
\end{equation}
Let $u=0$. Then $F=0$ and $\pa F=0$, so \eqref{eq:chain} and \eqref{eq:LF} show that $L_{G_1}(v)=L_F(v)=|\zeta|^2+|B\bar\xi|^2\ge0$, and that $L_{G_p}(v)=0$ for $p\ge2$.

Let $u\ne0$. Dividing $\pa F(v)$ by $|u|$, \eqref{eq:chain} reads $L_{G_p}(v)=pF^{p-1}\mathcal Q(\xi,\zeta)$ with
$$\mathcal Q(\xi,\zeta)=|\zeta|^2+|B\bar\xi|^2-2\Re\sum_i\bar u_iL_{f_i}(z;\xi)+(p-1)\big|\langle\zeta,e\rangle-c\big|^2,\qquad c=\langle e,B\bar\xi\rangle.$$
So $G_p$ has positive semidefinite Levi form at the point if and only if $\mathcal Q\ge0$ for all $\xi,\eta$. For fixed $\xi$, $\zeta$ is arbitrary with $\eta$. Write $\zeta=se+\zeta'$ with $\zeta'\perp e$. Then $\mathcal Q$ is smallest for $\zeta'=0$, and we must minimise $\psi(s)=|s|^2+(p-1)|s-c|^2$. Completing the square,
$$\psi(s)=p\,\Big|s-\fr{p-1}pc\Big|^2+\Big(p-1-\fr{(p-1)^2}p\Big)|c|^2=p\,\Big|s-\fr{p-1}pc\Big|^2+\Big(1-\fr1p\Big)|c|^2 .$$
So $\min_s\psi=(1-\fr1p)|c|^2$. As $|c|=|\langle B\bar\xi,e\rangle|$, the minimum of $\mathcal Q$ over $\zeta$ is $\Lev_{f,p}(z,w;\xi)$.
\end{proof}

\begin{Theorem}\label{thmHull}
Let $K$ be polynomially convex and $p\ge1$. If $\Lev_{f,p}(z,w;\cdot)\ge0$ at every point $(z,w)\in\widehat{\Gam}\setminus\Gam$, then $\Gam$ is polynomially convex.
\end{Theorem}
\begin{proof}
Recall $\widehat{\Gam}\subset K\times\C^N$. For $\ve>0$ let $\varphi_\ve=G_p+\ve(|z|^2+|w|^2)$, a $C^2$ function near $K\times\C^N$. By Proposition \ref{prop:levi} its Levi form is $\ge\ve|v|^2$ on the compact set $\widehat{\Gam}$, so by continuity $\varphi_\ve$ is plurisubharmonic on a neighbourhood $V_\ve$ of $\widehat\Gam$. Choose a Runge domain $D$, $\widehat\Gam\subset D\subset V_\ve$. By Theorem \ref{thm:psh}, for $(z,w)\in\widehat\Gam$,
$$|w-f(z)|^{2p}\le\varphi_\ve(z,w)\le\sup_{\Gam}\varphi_\ve\le\ve\sup_{\Gam}(|z|^2+|w|^2).$$
Let $\ve\to0$.
\end{proof}

Since the hull is usually unknown, one uses an outer estimate.

\begin{Corollary}\label{cor:loc}
(a) $\widehat\Gam\subset K\times\widehat{f(K)}$.

(b) If $\psi$ is plurisubharmonic on a Runge domain containing $\widehat\Gam$ and $\psi\le0$ on $\Gam$, then $\psi\le0$ on $\widehat\Gam$. In particular, if $g_i$ is holomorphic near $K$ and $|f_i|^2\le\Re(f_ig_i)$ on $K$, then $|w_i|\le|g_i(z)|$ on $\widehat\Gam$ (\cite{Chi}): use $\psi=|w_i|^2-\Re(w_ig_i(z))$.

Thus in Theorem \ref{thmHull} it suffices to have $\Lev_{f,p}\ge0$ on any compact set containing $\widehat\Gam$ described by (a), (b).
\end{Corollary}
\begin{proof}
(a) The projections of $\widehat\Gam$ lie in the hulls of $K$ and of $f(K)$. (b) The first claim is Theorem \ref{thm:psh}. For the second, $\psi$ is plurisubharmonic on $\Omega\times\C^N$ for a Runge domain $\Omega\supset K$ on which $g_i$ is holomorphic, and $\psi\le0$ on $\Gam$. On $\widehat\Gam$ we get $|w_i|^2\le\Re(w_ig_i)\le|w_i||g_i|$.
\end{proof}

\begin{Remark}\label{rem:gauge}
(1) Proposition \ref{prop:levi} and Theorem \ref{thmHull} hold for real $p\ge1$. If $N=1$ and the strict inequality $2|f_{\bar z}|^2>2\Re((\bar w-\bar f)f_{z\bar z})$ holds on a compact set containing $\widehat\Gam$, then, by compactness, \eqref{eq:one} holds for all large $p$.

(2) By Lemma \ref{lem:change}, $\Gam$ is polynomially convex if and only if the graph of $\tilde f=\Lambda f+b$ is, for $\Lambda,b$ holomorphic near $K$ with $\Lambda$ invertible. So condition (iii) may be tested on $\tilde f$. For constant $\Lambda$ this amounts to replacing $|w-f|^2$ by $|\Lambda(w-f)|^2$.
\end{Remark}

In the following situation no differential inequality is needed.

\begin{Proposition}\label{prop:twist}
Let $K$ be polynomially convex, $\Omega$ a Runge domain containing $K$, $\Lambda\colon\Omega\to GL_N(\C)$ and $b,k\colon\Omega\to\C^N$ holomorphic, and $f=\Lambda\,\bar k+b$. Then $\Gam$ is polynomially convex.
\end{Proposition}
\begin{proof}
Let $\Phi(z,w)=(z,\Lambda(z)^{-1}(w-b(z)))$. It is of the form in Lemma \ref{lem:change}, and it maps $\Gam$ onto the graph of $\bar k$ over $K$. So we may take $f=\bar k$. Expanding,
$$F(z,w):=|w-\ov{k(z)}|^2=\sum_i\big(|w_i|^2+|k_i(z)|^2-2\Re(w_ik_i(z))\big).$$
Here $|w_i|^2$ and $|k_i|^2$ are plurisubharmonic and $\Re(w_ik_i)$ is pluriharmonic, so $F$ is plurisubharmonic on $D=\Omega\times\C^N$. As $\widehat\Gam\subset K\times\C^N\subset D$ and $F=0$ on $\Gam$, Theorem \ref{thm:psh} gives $F\le0$ on $\widehat\Gam$. That is, $w=\ov{k(z)}$ on $\widehat\Gam$.
\end{proof}

The next result needs no smoothness. Its point is that the sweeping argument and the invertibility hypothesis in \cite[Thm.~3.8]{Pre} are not needed.

\begin{Proposition}\label{prop:peak}
Let $K\subset\C^N$ be polynomially convex, $f\colon K\to\C^M$ continuous. Suppose that for every $a\in K$ there is a polynomial $Q_a$ on $\C^N\times\C^M$ with $Q_a(a,w)=0$ for all $w$ and $\Re Q_a>0$ on $\Gam\setminus\{(a,f(a))\}$. Then $\Gam$ is polynomially convex. This applies with
\begin{equation}\label{eq:Qa}
Q_a(z,w)=\sum_j\beta_j(z_j-a_j)(w_j-f_j(a))+\sum_{i,j}\gamma_{ij}(z_i-a_i)(z_j-a_j)
\end{equation}
for constants depending on $a$; in particular, for $M=N$, $\beta_j=1$, $\gamma_{ij}=0$, under the monotonicity condition
\begin{equation}\label{eq:mon}
\Re\sum_j(z_j-a_j)\big(f_j(z)-f_j(a)\big)>0\qquad(z\ne a\text{ in }K).
\end{equation}
\end{Proposition}
\begin{proof}
Let $x=(a,\beta)\in\widehat\Gam$; then $a\in K$. Let $\mu$ be a representing probability measure for $x$ on $\Gam$. Since $e^{-Q_a}\in P(\Gam)$,
$$1=e^{-Q_a(x)}=\int\Re e^{-Q_a}\,d\mu\le\int e^{-\Re Q_a}\,d\mu\le1,$$
and $e^{-\Re Q_a}<1$ off $(a,f(a))$. Hence $\mu=\delta_{(a,f(a))}$ and $x=(a,f(a))$.
\end{proof}

For $M=N$ this removes the hypothesis ``$F_{\bar z}(a)$ invertible'' of \cite[Thm.~3.8]{Pre}, and it needs no regularity of $f$. For instance, it gives at once \cite[Lemma~5.3]{HW} and \cite[Thm.~3.6]{Pre} for polynomially convex $K$: if $f=\bar z+R$ with $|R(z)-R(z')|\le\kappa|z-z'|$, $\kappa<1$, then \eqref{eq:mon} holds. For arbitrary compact $K$ see Proposition \ref{prop:peakhull} and Theorem \ref{thmLip}.

The next result allows the holomorphic factors to degenerate. It extends the invertible case, Proposition \ref{prop:twist}, which is $h_1=\dots=h_N=1$.

\begin{Proposition}[Division]\label{prop:div}
Let $K$ be polynomially convex, $\Omega\supset K$ a Runge domain, and $\Lambda\colon\Omega\to GL_N(\C)$, $b,k,h\colon\Omega\to\C^N$ holomorphic. Let $D=\operatorname{diag}(h_1,\dots,h_N)$ and $f=\Lambda D\bar k+b$. Assume that for every $j$
\begin{equation}\label{eq:div}
\Re h_j\ge0\ \text{ on }K,\qquad\Re h_j>0\ \text{ on }K\setminus\{k_j=0\}.
\end{equation}
Then $\Gam$ is polynomially convex.
\end{Proposition}
\begin{proof}
By Lemma \ref{lem:change} we may assume $\Lambda=I$, $b=0$, so $f_j=\bar k_jh_j$. Let $x=(\alpha,\beta)\in\widehat\Gam$, so $\alpha\in K$, and let $\mu$ be a Jensen measure for $x$ on $\Gam$, which we regard as a measure on $K$. It represents $x$ on all functions holomorphic near $\widehat\Gam$, by Oka--Weil. Let $J=\{j:k_j(\alpha)\ne0\}$.

\emph{Step 1: $j\notin J$.} Let $Q=k_j(z)w_j$. Then $Q(x)=0$, and on $\Gam$ we have $\Re Q=|k_j|^2\Re h_j\ge0$. As in Proposition \ref{prop:peak}, $\int e^{-Q}\,d\mu=1$ forces $\Re Q=0$ $\mu$-almost everywhere. By \eqref{eq:div}, $\mu$ is carried by $\{k_j=0\}$. There $f_j=0$, so $\beta_j=\int f_j\,d\mu=0=f_j(\alpha)$.

\emph{Step 2: $j\in J$.} By \eqref{eq:div}, $h_j(\alpha)\ne0$. The Jensen inequality $\log|h_j(\alpha)|\le\int\log|h_j|\,d\mu$ gives $\mu(\{h_j=0\})=0$. For $t>0$ the function $\phi_{t,j}=1/(h_j+t)$ is holomorphic near $K$, since $\Re(h_j+t)\ge t$ on $K$, and $|h_j\phi_{t,j}|\le1$ on $K$.

\emph{Step 3.} Let $\psi_{t,j}=\phi_{t,j}$ for $j\in J$, $\psi_{t,j}=1$ for $j\notin J$, and let $q$ be a polynomial on $\C^N\times\C^N$. The function $q\big(z,(\psi_{t,j}(z)w_j)_j\big)$ is holomorphic near $\widehat\Gam$, so
$$q\big(\alpha,(\psi_{t,j}(\alpha)\beta_j)_j\big)=\int q\big(z,(\psi_{t,j}f_j)_j\big)\,d\mu.$$
For $j\in J$, $\psi_{t,j}f_j=\bar k_jh_j\phi_{t,j}$ is bounded by $|k_j|$ and tends to $\bar k_j$ off $\{h_j=0\}$ as $t\to0$. For $j\notin J$, $\psi_{t,j}f_j=f_j=0=\bar k_j$ $\mu$-almost everywhere, by Step 1. By dominated convergence,
$$q(\alpha,\beta')=\int q(z,\bar k)\,d\mu,\qquad \beta'_j=\beta_j/h_j(\alpha)\ (j\in J),\quad\beta'_j=0\ (j\notin J).$$
Thus $(\alpha,\beta')$ lies in the hull of the graph of $\bar k$ over $K$, which is polynomially convex by Proposition \ref{prop:twist}. So $\beta'=\ov{k(\alpha)}$, that is, $\beta_j=\ov{k_j(\alpha)}h_j(\alpha)=f_j(\alpha)$ for $j\in J$.
\end{proof}

\begin{Proposition}\label{prop:D}
Let $\mathcal C\subset\C^N$ be compact and convex, $f\in C^1$ near $\mathcal C$, and assume
\begin{equation}\label{eq:D}
\Re\big(v^{T}B_f(z)\bar v\big)>\big|v^{T}A_f(z)v\big|\qquad(z\in\mathcal C,\ v\in\C^N\setminus\{0\}).
\end{equation}
Then $B_f$ is invertible on $\mathcal C$, and $P[z,f;K]=C(K)$ for every compact $K\subset\mathcal C$. This holds if $f=\bar z+R$ with $\|dR\|\le\kappa<1$ on $\mathcal C$ (\cite{HW}).
\end{Proposition}
\begin{proof}
\emph{Invertibility.} If $B_f(z)\bar v=0$ with $v\ne0$, then \eqref{eq:D} reads $0>|v^TA_f(z)v|$, which is absurd. So $B_f$ is invertible on $\mathcal C$, as $v\mapsto\bar v$ is onto.

\emph{Polynomial convexity.} Let $a\ne z$ in $\mathcal C$, $v=z-a$ and $\gamma(t)=a+tv\in\mathcal C$. Since $\frac d{dt}f(\gamma(t))=A_f(\gamma)v+B_f(\gamma)\bar v$,
$$\sum_jv_j(f_j(z)-f_j(a))=\int_0^1\big(v^TA_f(\gamma)v+v^TB_f(\gamma)\bar v\big)\,dt.$$
Its real part is at least $\int_0^1\big(\Re(v^TB_f(\gamma)\bar v)-|v^TA_f(\gamma)v|\big)dt>0$. So \eqref{eq:mon} holds on the polynomially convex set $\mathcal C$, and $\Gam(\mathcal C)$ is polynomially convex by Proposition \ref{prop:peak}.

\emph{Density.} Since $B_f$ is invertible on $\mathcal C$, hence on a neighbourhood $U$ of $\mathcal C$, Lemma \ref{lem:tr} with $S=U$ shows that $\Gam(U)$ is a totally real $C^1$ manifold containing $\Gam(\mathcal C)$. Theorem \ref{thm:Iz} with $E=\emptyset$ gives $P(\Gam(\mathcal C))=C(\Gam(\mathcal C))$, so $P[z,f;\mathcal C]=C(\mathcal C)$ by the homeomorphism of $\mathcal C$ onto its graph. If $\phi\in C(K)$, Tietze's theorem extends it to $\mathcal C$, and polynomials in $z,f$ approximating the extension on $\mathcal C$ approximate $\phi$ on $K$.

\emph{The case $f=\bar z+R$.} Here $A_f=R_z$ and $B_f=I+R_{\bar z}$. Fix $v\ne0$, and choose $\theta$ such that the numbers $e^{i\theta}v^TR_zv$ and $e^{-i\theta}v^TR_{\bar z}\bar v$ have the same argument. With $dR(y)=R_zy+R_{\bar z}\bar y$ and $y=e^{i\theta}v$,
$$|v^TR_zv|+|v^TR_{\bar z}\bar v|=|v^TdR(e^{i\theta}v)|\le|v|\,\|dR\|\,|v|\le\kappa|v|^2.$$
Hence $\Re(v^TB_f\bar v)\ge|v|^2-|v^TR_{\bar z}\bar v|\ge(1-\kappa)|v|^2+|v^TR_zv|>|v^TA_fv|$.
\end{proof}
For $N=1$ condition \eqref{eq:D} says $\Re f_{\bar z}>|f_z|$; it is stronger than Duval's condition, but holds in all dimensions.

\begin{Remark}\label{rem:shift}
Let $U\in GL_N(\C)$ and let $S$ be a constant symmetric matrix. Proposition \ref{prop:D} and Theorem \ref{thmMon} hold, with the same proof, under the condition
\begin{equation}\label{eq:Dshift}
\Re\big(v^{T}UB_f(z)\,\bar v\big)>\big|v^{T}(UA_f(z)-S)\,v\big|\qquad(z\in\mathcal C,\ v\ne0):
\end{equation}
use $Q_a(z,w)=(z-a)^{T}U(w-f(a))-(z-a)^{T}S(z-a)$ in Proposition \ref{prop:peak}. So a constant part of $A_f$ does not matter. For $N=1$, $U=1$, $S=1$ and $f=\bar z+z$ this condition holds, while $|f_{\bar z}|=|f_z|$, so Duval's condition fails.
\end{Remark}

Now the singular set $Z=\{a\in K:\det B_f(a)=0\}$.

\begin{Proposition}[Stratification]\label{prop:strat}
Let $K$ be polynomially convex with $\Gam$ polynomially convex, $f\in C^1$. Let $K=A_0\supset A_1\supset\dots\supset A_r$ be polynomially convex compact sets, $A_r$ countable, such that each point of $\Gam(A_{i}\setminus A_{i+1})$ has a neighbourhood in which $\Gam(A_i)$ lies in a $C^1$ totally real submanifold. Then $P(\Gam)=C(\Gam)$.
\end{Proposition}
\begin{proof}
By (2.1) all $\Gam(A_i)$ are polynomially convex. Let $Y=\Gam(A_r)$, a countable compact set. For distinct $x,y\in Y$ the linear functionals $\ell$ on $\C^{2N}$ with $\ell(x)=\ell(y)$ form a hyperplane. There are countably many pairs, so by Baire's theorem some $\ell$ is injective on $Y$. Then $\ell(Y)\subset\C$ is compact and countable, so it has empty interior and connected complement, and $P(\ell(Y))=C(\ell(Y))$ by Lavrentiev's theorem. Since $\ell$ is a homeomorphism of $Y$ onto $\ell(Y)$, polynomials in $\ell$ are dense in $C(Y)$, so $P(Y)=C(Y)$.

Now we argue by descending induction on $i$. Suppose $P(\Gam(A_{i+1}))=C(\Gam(A_{i+1}))$. Apply Theorem \ref{thm:Iz} with $X=\Gam(A_i)$ and $E=\Gam(A_{i+1})$: the set $X$ is polynomially convex, and $X\setminus E=\Gam(A_i\setminus A_{i+1})$ satisfies the hypothesis. So $P(\Gam(A_i))=C(\Gam(A_i))$. For $i=0$ this is the claim.
\end{proof}

\begin{Corollary}\label{cor:strat}
Let $\Gam$ be polynomially convex. Then $P(\Gam)=C(\Gam)$ in each of the following cases.

(a) $P(\Gam(Z))=C(\Gam(Z))$; for instance, $Z$ countable.

(b) \emph{(Coordinate strata.)} For every $a\in K$, the columns of $B_f(a)$ with index in $\{i:a_i\neq0\}$ are linearly independent. \hfill{\rm (T)}

(c) $f=\Lambda D\bar k+b$ as in Proposition \ref{prop:div}, every point of $K$ is isolated in its $k$-fibre, and $\{h_j=0\}\subset\{k_j=0\}$ in $\Omega$ for each $j$. (The case $D=I$ is included.)
\end{Corollary}
\begin{proof}
(a) Let $U$ be a neighbourhood of $K$ on which $f$ is $C^1$, and $U'=U\cap\{\det B_f\ne0\}$, an open set. By Lemma \ref{lem:tr} with $S=U'$, $\Gam(U')$ is a totally real $C^1$ manifold, and it contains $\Gam\setminus\Gam(Z)$. Theorem \ref{thm:Iz} with $X=\Gam$ and $E=\Gam(Z)$ gives the claim. If $Z$ is countable, then $P(\Gam(Z))=C(\Gam(Z))$ by the first paragraph of the proof of Proposition \ref{prop:strat}.

(b) For $0\le d\le N$ let $A_d$ be the set of $z\in K$ with at least $d$ vanishing coordinates. A point has at least $d$ vanishing coordinates if and only if every product of $N-d+1$ distinct coordinates vanishes at it. So $A_d$ is the intersection of $K$ with the common zero set of finitely many polynomials, hence polynomially convex. Moreover $A_0=K$ and $A_N\subset\{0\}$. Let $a\in A_d\setminus A_{d+1}$ and $J=\{i:a_i=0\}$, so $|J|=d$. Near $a$ the coordinates with index outside $J$ do not vanish, so a point of $A_d$ near $a$ has exactly the vanishing coordinates $J$. Thus $A_d$ lies near $a$ in the coordinate plane $S=\{z_i=0,\ i\in J\}$, and $T_aS$ is spanned by the $e_i$, $i\notin J$. For $v\in T_aS\setminus\{0\}$ the vector $B_f(a)\bar v=\sum_{i\notin J}\bar v_i\,B_f(a)e_i$ is non-zero by (T). By Lemma \ref{lem:tr} the graph over $S$ is totally real at $(a,f(a))$, hence near this point, as this is an open condition. Proposition \ref{prop:strat} applies with $r=N$.

(c) Here $B_f=\Lambda D\,\ov{k'}$, so $B_f\bar v=\Lambda D\,\ov{k'v}$. The set of non-isolated points of fibres of $k$ is analytic in $\Omega$ by Cartan--Remmert, and it does not meet $K$. Replacing $\Omega$ by a smaller Runge domain, we may assume that $k$ has discrete fibres on $\Omega$.

Define analytic subsets $V_0=\Omega\supset V_1\supset\cdots$ of $\Omega$ as follows. Given $V_i$, let $V_{i+1}$ be the union of: the singular set of $V_i$; the sets $N_W$ defined below; and the sets $W\cap\{h_j=0\}$, for all irreducible components $W$ of $V_i$ and all $j$ with $h_j\not\equiv0$ on $W$. Here $N_W$ is the union of $\operatorname{Sing}W$ and the set of regular points of $W$ at which $k|_W$ is not an immersion. It is analytic. Indeed, near any point of $W$ choose generators $g_1,\dots,g_s$ of the ideal sheaf of $W$ (Oka's coherence theorem). At a regular point $z$ we have $T_zW=\ker dg(z)$, and $k|_W$ is an immersion at $z$ if and only if $\ker dg(z)\cap\ker dk(z)=0$, that is, the $(s+N)\times N$ matrix with rows $dg_1(z),\dots,dg_s(z),dk_1(z),\dots,dk_N(z)$ has rank $N$. So, locally,
$$N_W=\operatorname{Sing}W\cup\Big(W\cap\Big\{\operatorname{rank}\begin{pmatrix}dg\\ dk\end{pmatrix}<N\Big\}\Big),$$
and both sets on the right are analytic. The family of components is locally finite, so $V_{i+1}$ is analytic. Since $k|_W$ has discrete fibres, its generic rank on $W$ equals $\dim W$ by the rank theorem. So $N_W$ and the sets $W\cap\{h_j=0\}$ are proper analytic subsets of $W$, of dimension $<\dim W$, and $\operatorname{Sing}V_i$ has dimension $<\dim V_i$. Hence $\dim V_{i+1}<\dim V_i$, and $V_r$ is discrete for some $r\le N$. (Components of $V_i$ of dimension less than $\dim V_i$ need not be contained in $V_{i+1}$. This does no harm: only the property of the points of $V_i\setminus V_{i+1}$ proved below is used.) Each $V_i$ is the common zero set of a family of functions holomorphic on the Stein domain $\Omega$ (Cartan's Theorem A). These functions belong to $P(K)$ by Oka--Weil, so $A_i=K\cap V_i$ is polynomially convex; and $A_r$ is finite.

Let $a\in A_i\setminus A_{i+1}$. Near $a$, $S=V_i\setminus V_{i+1}$ is a complex manifold contained in one irreducible component $W$ of $V_i$, $k|_S$ is an immersion, and for each $j$ either $h_j\equiv0$ on $W$ or $h_j(a)\ne0$. Let $v\in T_aS\setminus\{0\}$ and $y=k'(a)v\ne0$, and choose $l$ with $y_l\ne0$. If $h_l(a)=0$, then $h_l\equiv0$ on $S$ near $a$, so $k_l\equiv0$ on $S$ near $a$ by hypothesis, and $y_l=0$, a contradiction. Hence $h_l(a)\ne0$, and $B_f(a)\bar v=\Lambda D\bar y\ne0$. Lemma \ref{lem:tr} shows that the graph over $S$ is totally real at $(a,f(a))$, and this persists near $a$ on $S$. Apply Proposition \ref{prop:strat}.
\end{proof}

\begin{Remark}\label{rem:fibre}
Condition (c) cannot be dropped. For $k=(z_1,z_1z_2)$ and $f=\bar k$ on the closed bidisk, $\Gam$ is polynomially convex by Proposition \ref{prop:twist}, $f=0$ on $\{z_1=0\}$, so $P[z,f;\ov\Delta^2]$ restricted to this disk is the disk algebra.
\end{Remark}

\section{Main theorems}

In this section $K\subset\C^N$ is compact, $m\in\N^N$, $c_j\in\N$ with $\gcd(c_j,m_j)=1$, and condition (ii) means
$$\text{(ii)}\quad z_j^{c_j}f_j(z)\in\Sect_{\pi/m_j}\qquad(m_j\ge2,\ z\in K,\ z_j\ne0).$$
All the proofs below have the same structure: one first proves $P[z,f;K]=C(K)$, using only properties of the graph $\Gam$, and then applies Theorem \ref{thmA}. We recall the first observation of the argument once. If $P(\Gam)=C(\Gam)$ and $\phi\in C(K)$, let $\tilde\phi(z,f(z))=\phi(z)$. Then $\tilde\phi\in C(\Gam)$, as $z\mapsto(z,f(z))$ is a homeomorphism of $K$ onto $\Gam$. A polynomial $p$ with $|p-\tilde\phi|<\ve$ on $\Gam$ satisfies $|p(z,f(z))-\phi(z)|<\ve$ on $K$. So
\begin{equation}\label{eq:GamK}
P(\Gam)=C(\Gam)\ \Longrightarrow\ P[z_1,\dots,z_N,f_1,\dots,f_N;K]=C(K).
\end{equation}
The application of Theorem \ref{thmA} is always with $M=N$ and $h_j(s,w)=w_j$, so that $\Phi_j(z)=z_j^{c_j}f_j(z)$, and condition (ii) is exactly \eqref{eq:sector}.

\begin{Theorem}\label{thmHess}
Let $K$ be polynomially convex and $f\in C^2$ near $K$ with values in $\C^N$. Assume (ii), one of the conditions (a), (b) of Corollary \ref{cor:strat} on the singular set, and
\begin{enumerate}
\item[(iii)] for some integer $p\ge1$, $\Lev_{f,p}(z,w;\cdot)\ge0$ at all points of $\widehat{\Gam}\setminus\Gam$; by Corollary \ref{cor:loc} it suffices that this hold for all $z\in K$, $w\in\widehat{f(K)}$ with $w\ne f(z)$, or on the smaller sets given there.
\end{enumerate}
Then $\Gam$ is polynomially convex and $P[z_1^{m_1},\dots,z_N^{m_N},f_1,\dots,f_N;K]=C(K)$.
\end{Theorem}
\begin{proof}
\emph{Step 1: $\Gam$ is polynomially convex.} Let $(z,w)\in\widehat\Gam\setminus\Gam$. By Corollary \ref{cor:loc}(a), $z\in K$ and $w\in\widehat{f(K)}$, and $w\ne f(z)$, because $z\in K$ and $(z,w)\notin\Gam$. So the two forms of (iii) both give $\Lev_{f,p}(z,w;\cdot)\ge0$ on $\widehat\Gam\setminus\Gam$, and Theorem \ref{thmHull} applies.

\emph{Step 2: $P(\Gam)=C(\Gam)$.} By Step 1 the graph is polynomially convex, and $f\in C^1$. So Corollary \ref{cor:strat} applies, under the assumption (a) or (b).

\emph{Step 3.} By \eqref{eq:GamK}, $P[z,f;K]=C(K)$. By (ii), Theorem \ref{thmA} gives the conclusion.
\end{proof}

\begin{Theorem}[Monotone case; replaces the use of Duval's theorem]\label{thmMon}
Let $\mathcal C$ be compact convex, $f\in C^1$ near $\mathcal C$ with \eqref{eq:D} on $\mathcal C$, and $K\subset\mathcal C$ compact with (ii). Then $P[z_1^{m_1},\dots,z_N^{m_N},f_1,\dots,f_N;K]=C(K)$.
\end{Theorem}
\begin{proof}
Proposition \ref{prop:D} gives $P[z,f;K]=C(K)$ for the compact subset $K$ of $\mathcal C$. Note that $K$ need not be polynomially convex: the polynomial convexity used is that of the graph over the convex set $\mathcal C$. By (ii), Theorem \ref{thmA} gives the conclusion; this theorem needs no assumption on $K$.
\end{proof}

\begin{Theorem}[Twisted antiholomorphic case, with division]\label{thmTwist}
Let $K$ be polynomially convex, $f=\Lambda D\bar k+b$ as in Proposition \ref{prop:div} (for instance $D=I$), with every point of $K$ isolated in its $k$-fibre and $\{h_j=0\}\subset\{k_j=0\}$ in $\Omega$, and assume (ii). Then $\Gam$ is polynomially convex and $P[z_1^{m_1},\dots,z_N^{m_N},f_1,\dots,f_N;K]=C(K)$.
\end{Theorem}
\begin{proof}
\emph{Step 1.} Proposition \ref{prop:div} shows that $\Gam$ is polynomially convex. In the case $D=I$ the conditions \eqref{eq:div} hold trivially, and this step is Proposition \ref{prop:twist}.

\emph{Step 2.} The hypotheses of Corollary \ref{cor:strat}(c) are part of the assumptions, so $P(\Gam)=C(\Gam)$.

\emph{Step 3.} By \eqref{eq:GamK}, $P[z,f;K]=C(K)$, and by (ii), Theorem \ref{thmA} gives the conclusion.
\end{proof}

\begin{Proposition}[Products]\label{prop:prod}
Let $K=K_1\times\dots\times K_N$, $K_j\subset\ov\Delta$ compact, $f_j=f_j(z_j)$ of class $C^1$ near $\ov\Delta$ with $|\pa f_j/\pa\bar z|>|\pa f_j/\pa z|$ on $\ov\Delta$, and $z_jf_j(z_j)\in\Sect_{\pi/m_j}$ on $K_j\setminus\{0\}$. Then $P[z^{[m]},f;K]=C(K)$.
\end{Proposition}
\begin{proof}
\emph{Step 1: one factor.} The function $f_j$ is $C^1$ near $\ov\Delta$ and $|\pa f_j/\pa\bar z|>|\pa f_j/\pa z|$ on $\ov\Delta$. These are exactly the hypotheses of Duval's theorem \cite[Th\'eor\`eme]{Du}. Its corollary \cite[Corollaire]{Du} gives $P[z_j,f_j;\ov\Delta]=C(\ov\Delta)$. Let $g\in C(K_j)$. By Tietze's theorem $g$ extends to $\ov\Delta$, so $g$ is a uniform limit on $\ov\Delta$, hence on $K_j$, of polynomials in $z_j,f_j(z_j)$. Thus $P[z_j,f_j;K_j]=C(K_j)$.

\emph{Step 2: the product.} Let $\mathcal A=P[z,f;K]$ and let $\pi_j\colon K\to K_j$ be the projections. If $g\in C(K_j)$ and $p_l(z_j,f_j(z_j))\to g$ uniformly on $K_j$, then $p_l(z_j,f_j(z_j))\to g\circ\pi_j$ uniformly on $K$; so $g\circ\pi_j\in\mathcal A$. Let $\mathcal B\subset\mathcal A$ be the algebra generated by the functions $g\circ\pi_j$, $g\in C(K_j)$, $1\le j\le N$. It contains the constants, it is closed under conjugation, as $\overline{g\circ\pi_j}=\bar g\circ\pi_j$, and it separates the points of $K$, as the coordinate functions belong to it. By the Stone--Weierstrass theorem $\mathcal B$ is dense in $C(K)$, so $\mathcal A=C(K)$.

\emph{Step 3.} We have $z_jf_j(z_j)\in\Sect_{\pi/m_j}$ for $z_j\in K_j\setminus\{0\}$, which is (ii) with $c_j=1$. Theorem \ref{thmA} gives the conclusion.
\end{proof}

\begin{Theorem}[Monomials of degree $n$]\label{thmNu}
Let $n\ge2$ and let $K$, $f$ satisfy the assumptions of Theorem \ref{thmHess}, \ref{thmMon} or \ref{thmTwist}, with (ii) replaced by
\begin{equation}\label{eq:sum}
\sum_{j=1}^Nz_jf_j(z)\in\Sect_{\pi/n}\qquad(z\in K\setminus\{0\}).
\end{equation}
Then $\nu_n,f$ separate the points of $K$ and $P[\nu_n,f_1,\dots,f_N;K]=C(K)$.
\end{Theorem}
\begin{proof}
Condition (ii) was used in those proofs only in the last step. So in each case the same argument gives $P[z,f;K]=C(K)$. Let $q(\zeta,w)=\sum_j\zeta_jw_j$. Then $q(\omega\zeta,w)=\omega q(\zeta,w)$, so $c=1$ is allowed in Theorem \ref{thmB}, and $\Phi(z)=q(z,f(z))=\sum_jz_jf_j(z)$. By \eqref{eq:sum}, Theorem \ref{thmB} gives the conclusion.
\end{proof}

\section{Examples}

$\ov\Delta$ denotes the closed unit disk.

\begin{Example}[Coupled Minsker theorem]\label{ex:1}
Let $m_j,s_j\in\N$ be coprime, $j=1,2$, and $0\le c<\pi/\max(m_1,m_2)$. Then
$$P\big[z_1^{m_1},z_2^{m_2},\ \bar z_1^{s_1}e^{cz_2},\ \bar z_2^{s_2}e^{cz_1};\ \ov\Delta^2\big]=C(\ov\Delta^2).$$
We use Theorem \ref{thmTwist} with $\Omega=\C^2$, $\Lambda=\operatorname{diag}(e^{cz_2},e^{cz_1})$, $k=(z_1^{s_1},z_2^{s_2})$, $D=I$, $b=0$, $c_j=s_j$.

The bidisk is polynomially convex; $\Lambda$ is holomorphic and invertible; and $f=\Lambda\bar k$, as $\ov{k_1}=\bar z_1^{s_1}$, $\ov{k_2}=\bar z_2^{s_2}$. A fibre of $k$ has at most $s_1s_2$ points, so every point is isolated in its fibre. It remains to check (ii). For $z_1\ne0$,
$$z_1^{s_1}f_1(z)=|z_1|^{2s_1}e^{cz_2},\qquad \arg e^{cz_2}=c\,\Im z_2\in[-c,c],$$
and $c<\pi/m_1$, so $z_1^{s_1}f_1\in\Sect_{\pi/m_1}$. The same holds for $f_2$. As $\gcd(s_j,m_j)=1$, the theorem applies.

For $c=0$ this is the product of Minsker's theorem. We have $\det B_f=\det\Lambda\cdot\ov{\det k'}=e^{c(z_1+z_2)}s_1s_2\,\bar z_1^{s_1-1}\bar z_2^{s_2-1}$. So $Z=\ov\Delta^2\cap\{z_1^{s_1-1}z_2^{s_2-1}=0\}$ is not finite unless $s_1=s_2=1$.
\end{Example}

\begin{Example}[An infinite singular set; a folded $f$]\label{ex:2}
Let $K=\{|z_1|\le3/5,\ |z_2|\le3/5\}$, $k=(z_1(1+z_2),z_2(1+z_1))$, $f=\bar k$. Then
$$P\big[z_1^{4},z_2^{4},\ \bar z_1(1+\bar z_2),\ \bar z_2(1+\bar z_1);K\big]=C(K).$$
We use Theorem \ref{thmTwist} with $\Lambda=I$, $D=I$, $b=0$, $c_j=1$, $m_j=4$.

\emph{Condition (ii).} We have $z_1f_1=|z_1|^2(1+\bar z_2)$. If $|s|\le\rho<1$, then $1+s$ lies in the disk of centre $1$ and radius $\rho$, seen from $0$ under the half-angle $\arcsin\rho$; so $|\Arg(1+s)|\le\arcsin\rho$. With $\rho=3/5$, and since $3/5<\sqrt2/2=\sin(\pi/4)$, we get $|\Arg(1+\bar z_2)|<\pi/4$. The same holds for $z_2f_2=|z_2|^2(1+\bar z_1)$.

\emph{Fibres.} Let $k(z)=(A_1,A_2)$. As $k_1-k_2=z_1-z_2$, we have $z_1=z_2+\delta$ with $\delta=A_1-A_2$, and then $z_2(1+z_2+\delta)=A_2$. This quadratic equation has at most two solutions, so the fibres have at most two points. Theorem \ref{thmTwist} applies.

\emph{The singular set.} We have
$$k'=\begin{pmatrix}1+z_2&z_1\\ z_2&1+z_1\end{pmatrix},\qquad\det k'=(1+z_1)(1+z_2)-z_1z_2=1+z_1+z_2 .$$
So $Z=K\cap\{z_1+z_2=-1\}$, a two-dimensional convex set containing $(-\frac12,-\frac12)$. The map $k$ is not injective on $K$: $t+t^2=s+s^2$ for $s=-1-t$, so $k(t,t)=k(-1-t,-1-t)$ for $-3/5\le t\le-2/5$. Finally, $k|_Z$ is an immersion, as $k'(1,-1)^T=(1+z_2-z_1,\ z_2-z_1-1)$ has two entries differing by $2$.
\end{Example}

\begin{Example}[Monotone case]\label{ex:3}
Let $\mathcal C$ be the closed unit ball of $\C^2$, $K\subset\mathcal C$ compact, $m_1,m_2\le9$, and
$$f=\big(\bar z_1+\tfrac13z_1|z_2|^2,\ \bar z_2+\tfrac13z_2|z_1|^2\big).$$
Then $P[z_1^{m_1},z_2^{m_2},f_1,f_2;K]=C(K)$. We use Theorem \ref{thmMon} with $c_j=1$.

\emph{Condition \eqref{eq:D}.} Let $R=f-\bar z$, so $R_1=\frac13z_1z_2\bar z_2$. Then
$$dR_1(v)=\tfrac13\big(|z_2|^2v_1+z_1\bar z_2v_2+z_1z_2\bar v_2\big),\qquad |dR_1(v)|\le\tfrac13\big(|z_2|^2|v_1|+2|z_1z_2||v_2|\big),$$
and similarly $|dR_2(v)|\le\frac13(|z_1|^2|v_2|+2|z_1z_2||v_1|)$. On $\mathcal C$ we have $|z_j|^2\le1$ and $2|z_1z_2|\le|z|^2\le1$. So $|dR_j(v)|\le\frac13(|v_1|+|v_2|)\le\frac{\sqrt2}3|v|$, and
$$|dR(v)|^2\le2\cdot\tfrac29|v|^2=\tfrac49|v|^2 .$$
Thus $\|dR\|\le\frac23<1$ on $\mathcal C$, and \eqref{eq:D} holds by Proposition \ref{prop:D}.

\emph{Condition (ii).} For $z_1\ne0$, $z_1f_1=|z_1|^2+\frac13z_1^2|z_2|^2=|z_1|^2(1+s)$ with $s=\frac13(z_1/|z_1|)^2|z_2|^2$ and $|s|\le\frac13$. So its argument has modulus at most $\arcsin\frac13$, and $\frac13<\sin(\pi/9)\approx0.342$. Hence (ii) holds for $m_1\le9$, and likewise for $z_2f_2$.

Note that $f$ is not pluriharmonic and that $K$ need not be polynomially convex.
\end{Example}

\begin{Example}[The Hermitian criterion, with degeneration on coordinate axes]\label{ex:4}
Let $\mu\ge2$, $0<\alpha$, $0\le c$, $\alpha'=\alpha(1+c)<1$,
$$f_1=\bar z_1^{\mu}+\alpha\bar z_1z_1^{\mu-1}(1+cz_2^{\mu}),\qquad f_2=\bar z_2^{\mu}+\alpha\bar z_2z_2^{\mu-1}(1+cz_1^{\mu}),$$
and put $a=(1+\alpha')^2/(1-\alpha')$, $C_0=a+1+\alpha'$. Assume
\begin{equation}\label{eq:ex4}
(\mu-\alpha')^2\ge2C_0\big[\alpha'(\mu-1)+\alpha c\mu\big],
\end{equation}
and that $m_1,m_2$ are coprime to $\mu$ with $\arcsin\alpha'<\pi/m_j$. Then $P[z_1^{m_1},z_2^{m_2},f_1,f_2;\ov\Delta^2]=C(\ov\Delta^2)$. For example, $\mu=3$, $\alpha=1/10$, $c=1$, $m_1=2$, $m_2=5$: then $\alpha'=1/5$, $a=9/5$, $C_0=3$, and \eqref{eq:ex4} reads $7.84\ge4.2$.

Verification. Write $f_1=\bar z_1^\mu(1+\theta_1)$, $|\theta_1|=\alpha|1+cz_2^\mu|\le\alpha'$ for $z_1\ne0$. Then (ii) holds with $c_j=\mu$, since $z_1^\mu f_1=|z_1|^{2\mu}(1+\theta_1)$. With $g_j=az_j^\mu$ we have $|f_j|^2\le(1+\alpha')^2|z_j|^{2\mu}\le a(1-\alpha')|z_j|^{2\mu}\le\Re(f_jg_j)$, so by Corollary \ref{cor:loc} on $\widehat\Gam$ we have $|w_j|\le a|z_j|^\mu$ and $|u_j|\le C_0|z_j|^\mu$. Put $t_j=|z_j|^{\mu-1}$. Here $B_f=\operatorname{diag}(b_1,b_2)$ with $|b_j|\ge(\mu-\alpha')t_j$, $b_1=\mu\bar z_1^{\mu-1}+\alpha z_1^{\mu-1}(1+cz_2^\mu)$, and
$$L_{f_1}(z;\xi)=\bar\xi_1\big[\alpha(\mu-1)z_1^{\mu-2}(1+cz_2^\mu)\xi_1+\alpha c\mu z_1^{\mu-1}z_2^{\mu-1}\xi_2\big],$$
similarly for $f_2$. Since $|z_j|\le1$,
\begin{align*}
2\Big|\sum_i\bar u_iL_{f_i}\Big|&\le2C_0\alpha'(\mu-1)\big(t_1^2|\xi_1|^2+t_2^2|\xi_2|^2\big)+4C_0\alpha c\mu\,t_1t_2|\xi_1\xi_2|\\
&\le2C_0\big[\alpha'(\mu-1)+\alpha c\mu\big]\big(t_1^2|\xi_1|^2+t_2^2|\xi_2|^2\big),
\end{align*}
while $|B_f\bar\xi|^2\ge(\mu-\alpha')^2(t_1^2|\xi_1|^2+t_2^2|\xi_2|^2)$. So $\Lev_{f,1}\ge0$ on $\widehat\Gam$. Finally $\det B_f=b_1b_2$ vanishes exactly on $z_1z_2=0$, and (T) holds as $B_f$ is diagonal with $b_j\neq0$ when $z_j\ne0$. Theorem \ref{thmHess} applies.

This $f$ is not of the form treated in Theorem \ref{thmMon}: for $\mu=3$, $\omega=e^{2\pi i/3}$, $z=(\omega,0)$, $a=(\bar\omega,0)$, the left side of \eqref{eq:mon} equals $-3\alpha<0$.
\end{Example}

\begin{Example}[Monomials of degree $n$]\label{ex:5}
Let $\mathcal B$ be the closed unit ball in $\C^2$, $b,c\ge0$, $\arctan b+c<\pi/n$, $h=e^{cz_1}$, and
$$f=\big((\bar z_1+b\bar z_2)h,\ (\bar z_2-b\bar z_1)h\big).$$
Then $P[\nu_n,f_1,f_2;K]=C(K)$ for each polynomially convex compact $K\subset\mathcal B$. We use Theorem \ref{thmNu} in the setting of Theorem \ref{thmTwist}.

We have $f=\Lambda\bar k$ with $\Lambda=hI$ and $k=(z_1+bz_2,\ z_2-bz_1)$. Here $\Lambda$ is invertible on $\C^2$, and $k$ is linear with determinant $1+b^2\ne0$, so its fibres are points. Further
$$\sum_jz_jf_j=h\big(|z_1|^2+bz_1\bar z_2+|z_2|^2-b\bar z_1z_2\big)=h\big(|z|^2+2ib\,\Im(z_1\bar z_2)\big).$$
Since $2|\Im(z_1\bar z_2)|\le2|z_1||z_2|\le|z|^2$, the second factor has argument of modulus at most $\arctan b$ for $z\ne0$. Also $\arg h=c\,\Im z_1$ has modulus at most $c$ on $\mathcal B$. So \eqref{eq:sum} holds, as $\arctan b+c<\pi/n$. For $b>0$ and, say, $K=\mathcal B$, condition \eqref{eq:sector} with $h_j=w_j$ fails, since $f_1=0$ on $\bar z_1=-b\bar z_2$; so Theorem \ref{thmA} does not apply in this form.
\end{Example}

\begin{Example}[Division: $\det B_f$ vanishes and $\Lambda D$ is singular]\label{ex:6}
Let $\theta>0$, $c\ge0$, $m\ge2$ with $\theta+c<\pi/m$, let $T=\{re^{i\varphi}:0\le r\le1,\ |\varphi|\le\theta\}$ and $K=T\times T$. Then
$$P\big[z_1^{m},z_2^{m},\ |z_1|^2e^{cz_2},\ |z_2|^2e^{cz_1};\ K\big]=C(K).$$
Use Theorem \ref{thmTwist} with $k=z$, $\Lambda=I$, $b=0$, $h=(z_1e^{cz_2},z_2e^{cz_1})$. The set $T$ is compact with connected complement, so $K$ is polynomially convex; and $\{h_j=0\}=\{z_j=0\}=\{k_j=0\}$. On $K$ with $z_1\ne0$, $\arg h_1=\varphi_1+c\,\Im z_2$ has modulus at most $\theta+c<\pi/m\le\pi/2$, so \eqref{eq:div} holds, and (ii) holds with $c_j=1$ since $z_1f_1=|z_1|^2h_1$. Here $f$ is not of the form of Proposition \ref{prop:twist}, as $D$ is singular on $\{z_1z_2=0\}$.
\end{Example}

\begin{Example}[A cyclic quotient singularity]\label{ex:7}
Let $K\subset\C^2$ be polynomially convex. Then
$$P\big[z_1^3,\ z_1z_2,\ z_2^3,\ \bar z_1,\ \bar z_2^2;\ K\big]=C(K).$$
Use Theorem \ref{thmC} with $n=3$, $a=(1,2)$: the monomials $z_1^\alpha z_2^\gamma$ with $\alpha+2\gamma\equiv0\pmod3$ are generated by $z_1^3,z_1z_2,z_2^3$. Let $q=\zeta_1w_1+\zeta_2^2w_2$. Then $q(T_a\zeta,w)=\omega q$ since $\omega^4=\omega$, and $\Phi=|z_1|^2+|z_2|^4$ is positive on $K\setminus\{0\}$. Finally $P[z,\bar z_1,\bar z_2^2;K]=C(K)$ by Proposition \ref{prop:twist} and Corollary \ref{cor:strat}(c), with $k=(z_1,z_2^2)$.
\end{Example}

\section{The case $N=1$}\label{sec:one}

As a consequence of Theorem \ref{thmHess} with $N=1$, we obtain the following theorem, which parallels a result in \cite{CT}.
\begin{Theorem}\label{thmOne}
Let $K\subset\C$ be polynomially convex, $f\in C^2$ near $K$, $Z=\{f_{\bar z}=0\}\cap K$ finite or countable. Let $\gcd(c,n)=1$ and $z^cf(z)\in\Sect_{\pi/n}$ for $z\in K\setminus\{0\}$. If for some integer $p\ge1$
\begin{equation}\label{eq:one}
\Big(2-\fr1p\Big)|f_{\bar z}(z)|^2\ge2\Re\Big(\big(\bar w-\ov{f(z)}\big)f_{z\bar z}(z)\Big)\quad\text{for }(z,w)\in\widehat{\Gam},
\end{equation}
then $\Gam$ is polynomially convex and $P[z^n,f;K]=C(K)$.
\end{Theorem}
\begin{proof}
Theorem \ref{thmHull} with $N=1$, Corollary \ref{cor:strat}(a), and Theorem \ref{thmA} with $\Phi=z^cf$.
\end{proof}

\begin{Example}\label{ex:vd}
Let $0<\alpha<1$, $\mu\ge2$, $n\ge3$, $f=\bar z^\mu+\alpha\bar zz^{\mu-1}$, $a=(1+\alpha)^2/(1-\alpha)$, and assume
$$(2-1/p)(\mu-\alpha)^2\ge2\alpha(\mu-1)(a+1+\alpha)\quad\text{for some }p\in\N.$$
Then $P[z^n,f;K]=C(K)$ in each of the following cases:

(1) $K=\{re^{i\varphi}:0\le r\le1,\ \cos((\mu-1)\varphi)\ge\kappa\}$ with $\kappa=\fr{1-\alpha}{1+\alpha}\cot\fr\pi n$;

(2) $K=\ov\Delta$, $\gcd(\mu,n)=1$ and $\arcsin\alpha<\pi/n$.

\emph{The differential inequality.} We use Theorem \ref{thmOne}. For $z\ne0$ write $f=\bar z^\mu(1+\theta)$, $\theta=\alpha z^{\mu-1}/\bar z^{\mu-1}$, so $|\theta|=\alpha$. With $g=az^\mu$,
$$|f|^2\le(1+\alpha)^2|z|^{2\mu}=a(1-\alpha)|z|^{2\mu}\le a|z|^{2\mu}\Re(1+\theta)=\Re(fg).$$
By Corollary \ref{cor:loc}(b), $|w|\le a|z|^\mu$ on $\widehat\Gam$, so
$$|\bar w-\bar f|\le|w|+|f|\le(a+1+\alpha)|z|^\mu\quad\text{on }\widehat\Gam.$$
Further $f_{\bar z}=\mu\bar z^{\mu-1}+\alpha z^{\mu-1}$ and $f_{z\bar z}=\alpha(\mu-1)z^{\mu-2}$. So $|f_{\bar z}|\ge(\mu-\alpha)|z|^{\mu-1}$, which shows $Z=\{0\}$ (as $\alpha<1<\mu$), and on $\widehat\Gam$
$$2\Re\big((\bar w-\bar f)f_{z\bar z}\big)\le2\alpha(\mu-1)(a+1+\alpha)|z|^{2\mu-2}\le(2-\tfrac1p)(\mu-\alpha)^2|z|^{2\mu-2}\le(2-\tfrac1p)|f_{\bar z}|^2 ,$$
which is \eqref{eq:one}.

\emph{The angle.} Let $z=re^{i\varphi}\ne0$ and $\psi=(\mu-1)\varphi$. Then
$$zf=r^2\big(\bar z^{\mu-1}+\alpha z^{\mu-1}\big)=r^{\mu+1}\big((1+\alpha)\cos\psi-i(1-\alpha)\sin\psi\big).$$
In case (1), $\cos\psi\ge\kappa>0$ as $n\ge3$, so $\Re(zf)>0$, and since $|\sin\psi|<1$,
$$|\tan\arg(zf)|=\fr{(1-\alpha)|\sin\psi|}{(1+\alpha)\cos\psi}<\fr{1-\alpha}{(1+\alpha)\kappa}=\tan\fr\pi n .$$
So the condition of Theorem \ref{thmOne} holds with $c=1$. In case (2), $z^\mu f=|z|^{2\mu}(1+\theta')$ with $\theta'=\alpha(z/|z|)^{2\mu-2}$, so $|\arg(z^\mu f)|\le\arcsin\alpha<\pi/n$, and the condition holds with $c=\mu$.
\end{Example}

For $f=\bar zh$ with $h$ holomorphic no differential inequality is needed.

\begin{Proposition}\label{prop:division}
Let $K\subset\C$ be polynomially convex, $h$ holomorphic near $K$, $n\ge2$, and $h(K\setminus\{0\})\subset\Sect_{\pi/n}$. Then the graph of $\bar zh$ over $K$ is polynomially convex and $P[z^n,\bar zh;K]=C(K)$.
\end{Proposition}
\begin{proof}
Let $f=\bar zh$. If $h$ has no zero on $K$, this is Theorem \ref{thmTwist} with $\Lambda=h$, $k=z$, $D=I$, $b=0$. Otherwise $0\in K$ and $h(0)=0$, the only zero of $h$ on $K$. Then $\Re h\ge0$ on $K$ and $\Re h>0$ off $\{z=0\}$, so $\Gam$ is polynomially convex by Proposition \ref{prop:div} with $k=z$, $D=h$. As $Z=\{0\}$, Corollary \ref{cor:strat}(a) gives $P(\Gam)=C(\Gam)$, and Theorem \ref{thmA} with $c=1$ concludes.
\end{proof}
This removes conditions (b) and (c) of \cite[Cor.~3.5]{CT}; (c) in fact follows from (a), as $\Re(h(a)+h(-a))>0$.

The theorems of \cite{CT} based on Duval's theorem are unaffected. There, \cite[Corollaire]{Du} gives $P[z,f;\ov D]=C(\ov D)$ on a closed disk $\ov D$ (after an affine change of variable, which preserves the condition $|f_z|<|f_{\bar z}|$). By Tietze's theorem, as in Step 1 of the proof of Proposition \ref{prop:prod}, this gives $P[z,f;K]=C(K)$ for compact $K\subset\ov D$. The gluing argument is Lemma \ref{lem:rot}.

\section{Lipschitz perturbations of $\bar z$ on arbitrary compact sets}\label{sec:HW}

H\"ormander and Wermer \cite[Thm.~1.2]{HW} proved that $P[z_1,\dots,z_N,\bar z_1+R_1,\dots,\bar z_N+R_N;X]=C(X)$ for every compact $X\subset\C^N$, if $R$ is of class $C^{N+1}$ and satisfies a Lipschitz condition with constant $k<1$ on a neighbourhood of $X$. Weinstock \cite{Wei} gave a proof by Cauchy--Fantappi\`e kernels for $R\in C^1$ \cite[Thm.~14.4]{AW}. Following a remark of Browder, Cole and Wermer, the Lipschitz condition on $X$ alone suffices. We use the result in this form.

\begin{Theorem}[{\cite[Thm.~1.6.10, Cor.~1.6.12]{St}}]\label{thmLip}
Let $X\subset\C^N$ be compact and let $R\colon X\to\C^N$ satisfy
\begin{equation}\label{eq:Lip}
|R(z)-R(z')|\le k|z-z'|\qquad(z,z'\in X),\qquad k<1.
\end{equation}
Then $P[z_1,\dots,z_N,\bar z_1+R_1,\dots,\bar z_N+R_N;X]=C(X)$, and the graph of $\bar z+R$ over $X$ is polynomially convex.
\end{Theorem}

The set $X$ is arbitrary; it need not be polynomially convex. The proof in \cite{St} extends $R$ to $\C^N$ with the same constant (Kirszbraun's theorem \cite{Kir,Val}) and passes from smooth to Lipschitz maps by mollification. The polynomial convexity of the graph can also be seen directly, by the argument of \cite[Lemma~5.3]{HW}. We record that argument in an abstract form, since it is used again in Proposition \ref{prop:locDuval}.

\begin{Proposition}[Peak functions over the hull]\label{prop:peakhull}
Let $K\subset\C^N$ be compact (not necessarily polynomially convex), $f\colon K\to\C^M$ continuous, and $W\supset\widehat K$ open. Suppose that for every $a\in\widehat K$ there is a function $Q_a$ holomorphic on $W\times\C^M$ with $Q_a(a,w)=0$ for all $w\in\C^M$ and
$$\Re Q_a>0\quad\text{on }\Gam\setminus(\{a\}\times\C^M).$$
Then $\Gam$ is polynomially convex.
\end{Proposition}
\begin{proof}
Let $x=(a,\beta)\in\widehat\Gam$. The projection of $\widehat\Gam$ to $\C^N$ lies in $\widehat K$, so $a\in\widehat K$, and $\widehat\Gam$ is a polynomially convex compact subset of $W\times\C^M$. By the Oka--Weil theorem there are polynomials $p_l$ with $p_l\to e^{-Q_a}$ uniformly on $\widehat\Gam$. Let $\mu$ be a probability measure on $\Gam$ representing $x$ for $P(\Gam)$. Then
$$\int e^{-Q_a}\,d\mu=\lim_l\int p_l\,d\mu=\lim_lp_l(x)=e^{-Q_a(x)}=1 .$$
The set $\Gam\cap(\{a\}\times\C^M)$ is empty if $a\notin K$ and equals $\{(a,f(a))\}$ if $a\in K$; there $Q_a=0$. Hence $\Re Q_a\ge0$ on $\Gam$, so $|e^{-Q_a}|=e^{-\Re Q_a}\le1$ on $\Gam$, with equality only on $\Gam\cap(\{a\}\times\C^M)$. From
$$1=\Re\int e^{-Q_a}\,d\mu\le\int e^{-\Re Q_a}\,d\mu\le1$$
we get $e^{-\Re Q_a}=1$ $\mu$-almost everywhere, so $\mu$ is carried by $\Gam\cap(\{a\}\times\C^M)$. As $\mu$ is a probability measure, this set is not empty. Thus $a\in K$ and $\mu=\delta_{(a,f(a))}$, and $x=\int(z,w)\,d\mu=(a,f(a))\in\Gam$.
\end{proof}
For polynomially convex $K$ and polynomial $Q_a$ this is Proposition \ref{prop:peak}. The conclusion does not say that $K$ is polynomially convex: the points of $\widehat K\setminus K$ simply do not lie under $\widehat\Gam$. For $f=\bar z+R$ as in Theorem \ref{thmLip}, let $\tilde R$ be a $k$-Lipschitz extension of $R$ to $\C^N$, and for $a\in\widehat X$ put
$$Q_a(z,w)=\sum_j(z_j-a_j)\big(w_j-\bar a_j-\tilde R_j(a)\big).$$
For $z\in X$, the Cauchy--Schwarz inequality gives
$$\Re Q_a(z,f(z))=|z-a|^2+\Re\sum_j(z_j-a_j)\big(R_j(z)-\tilde R_j(a)\big)\ge(1-k)|z-a|^2 .$$
So Proposition \ref{prop:peakhull} gives the polynomial convexity in Theorem \ref{thmLip}.

\begin{Corollary}\label{corHW}
Let $K\subset\C^N$ be compact, not necessarily polynomially convex, and let $f=\Lambda(\bar z+H+R)$ on $K$. Here $\Lambda\in GL_N(\C)$ is constant, $H$ is holomorphic on a neighbourhood of $\widehat K$, and $R$ satisfies \eqref{eq:Lip} on $K$.
\begin{enumerate}
\item If (ii) holds, then $P[z_1^{m_1},\dots,z_N^{m_N},f_1,\dots,f_N;K]=C(K)$. When $\Lambda=I$ and $H=0$, (ii) holds with $c_j=1$ if $|R_j(z)|\le s_j|z_j|$ on $K$ with $s_j<\sin(\pi/m_j)$.
\item If \eqref{eq:sum} holds, then $P[\nu_n,f_1,\dots,f_N;K]=C(K)$. When $\Lambda=I$ and $H=0$, \eqref{eq:sum} holds if $|\sum_jz_jR_j(z)|\le s|z|^2$ on $K$ with $s<\sin(\pi/n)$.
\end{enumerate}
\end{Corollary}
\begin{proof}
By Oka--Weil, $H\in P(\widehat K)$. As $\Lambda$ is constant and invertible, $P[z,f;K]=P[z,\bar z+R;K]$, which is $C(K)$ by Theorem \ref{thmLip}. Theorems \ref{thmA} and \ref{thmB} need no assumption on $K$ and give the two claims. For the sufficient conditions: if $z_j\ne0$, then $z_jf_j=|z_j|^2(1+R_j/\bar z_j)$ with $|R_j/\bar z_j|\le s_j<1$. So $|\Arg(z_jf_j)|\le\arcsin s_j<\pi/m_j$. Likewise $\sum_jz_jf_j=|z|^2+\sum_jz_jR_j$ and $|\Arg\sum_jz_jf_j|\le\arcsin s<\pi/n$ for $z\ne0$.
\end{proof}
The gauge $\Lambda,H$ does not change $P[z,f;K]$, but it matters for the powers. $H$ need not belong to the algebra generated by $z^{[m]}$ or $\nu_n$, and $\Lambda$ mixes the components.

\begin{Example}[Inversion on a shell]\label{ex:shell}
Let $N\ge2$, $0<r_0\le1$, and let $K$ be any compact subset of the shell $\{r_0\le|z|\le1\}$, for instance the sphere $\{|z|=1\}$ or the shell itself; such $K$ are in general not polynomially convex. Let $\ve>0$, $\kappa=\ve/r_0^2<1$, and
$$f(z)=\bar z+\ve\,\fr{z}{|z|^2}.$$
Then $P[z_1^{m_1},\dots,z_N^{m_N},f;K]=C(K)$ whenever $\kappa<\sin(\pi/m_j)$ for all $j$ with $m_j\ge2$, and $P[\nu_n,f;K]=C(K)$ whenever $\kappa<\sin(\pi/n)$.

We apply Corollary \ref{corHW} with $\Lambda=I$, $H=0$, $R(z)=\ve z/|z|^2$. Expanding the squares one checks the identity
$$\Big|\fr z{|z|^2}-\fr{z'}{|z'|^2}\Big|=\fr{|z-z'|}{|z||z'|},$$
so \eqref{eq:Lip} holds on $K$ with $k=\ve/r_0^2=\kappa<1$, without any convexity of $K$. Further $|R_j(z)|=\ve|z_j|/|z|^2\le\kappa|z_j|$, and $|\sum_jz_jR_j|=\ve|\sum_jz_j^2|/|z|^2\le\ve\le\kappa|z|^2$ on $K$.

Theorem \ref{thmLip} alone gives no information on powers. Theorem \ref{thmMon} does not apply either. It needs a compact convex set $\mathcal C\supset K$ with $f\in C^1$ near $\mathcal C$, and for $K$ the sphere every such $\mathcal C$ contains $0$, where $f$ is unbounded. Theorems \ref{thmHess} and \ref{thmTwist} require $K$ polynomially convex, which the sphere is not.
\end{Example}

\begin{Remark}\label{rem:HWcomp}
(1) \emph{Comparison with Theorem \ref{thmMon}.} Suppose \eqref{eq:Lip} holds on a compact convex set $\mathcal C\supset K$ with non-empty interior, and $R\in C^1$ near $\mathcal C$. Then $\|dR\|\le k$ on $\mathcal C$, by $|dR(a)\theta|=\lim_{t\to0^+}|R(a+t\theta)-R(a)|/t$ at interior points and continuity. This implies \eqref{eq:D} (proof of Proposition \ref{prop:D}). So in that case Corollary \ref{corHW}(1) with $\Lambda=I$, $H=0$ is contained in Theorem \ref{thmMon}. The point of Corollary \ref{corHW} is that \eqref{eq:Lip} is needed only on $K$, without convexity, as in Example \ref{ex:shell}.

(2) \emph{$N=1$.} For the disk, Theorem \ref{thmLip} is Wermer's theorem \cite{W}, \cite[Thm.~14.3]{AW}. If $R\in C^1$ with $\|dR\|\le k$ near $\ov\Delta$, then $\|dR\|=|R_z|+|R_{\bar z}|$ gives $|f_z|=|R_z|<1-|R_{\bar z}|\le|f_{\bar z}|$. So the $C^1$ case on the disk is contained in Duval's theorem \cite{Du}. The Lipschitz case, and the case of arbitrary compact sets, are not.

(3) \emph{Sharpness.} The condition $k<1$ cannot be relaxed to $k=1$: take $R=-\bar z$ on $\ov\Delta$ \cite[p.~16]{HW}, \cite[p.~55]{St}.
\end{Remark}

The next statement is a direct consequence of \cite[Thm.~5.4]{HW}.

\begin{Corollary}[Stability]\label{cor:stab}
Let $K\subset\C^N$ be compact, $U\supset K$ open, $f\in C^{N+1}(U,\C^N)$ with $\det B_f\ne0$ on $K$, and $P[z,f;K]=C(K)$. Then $P[z,g;K]=C(K)$ for every $g\in C^{N+1}(U,\C^N)$ sufficiently close to $f$ in the $C^2$ topology of \cite[Thm.~5.4]{HW}. If such a $g$ satisfies (ii) (resp. \eqref{eq:sum}), then $P[z^{[m]},g;K]=C(K)$ (resp. $P[\nu_n,g;K]=C(K)$).
\end{Corollary}
\begin{proof}
Let $M\subset U$ be an open neighbourhood of $K$ on which $\det B_f\ne0$, a manifold of real dimension $2N$, and $\varphi=(z,f)\colon M\to\C^{2N}$. Then $\varphi\in C^{r'}$ with $r'=N+1=\max(2,1+\frac{2N}2)$. Hypothesis (i) of \cite[Thm.~5.4]{HW} is $P[z,f;K]=C(K)$. For (ii), use $dz,d\bar z$ as coordinates on the complexified tangent space: the complex linear extension of $d\varphi$ has matrix $\begin{pmatrix}I&0\\A_f&B_f\end{pmatrix}$, which is injective if and only if $\det B_f\ne0$. So \cite[Thm.~5.4]{HW} gives $P[z,g;K]=C(K)$ for $\psi=(z,g)$ close to $\varphi$, that is, for $g$ close to $f$. Theorems \ref{thmA} and \ref{thmB} conclude.
\end{proof}

\begin{Remark}[Instability when $Z\ne\emptyset$]\label{rem:instab}
The hypothesis $\det B_f\ne0$ cannot be dropped, as \cite[Ex.~6.2]{HW} shows for $N=1$. In two variables, Example \ref{ex:1} with $c=0$, $s_1=3$, $s_2=1$ and $\gcd(m_1,3)=1$ gives
$$P[z_1^{m_1},z_2^{m_2},\bar z_1^3,\bar z_2;\ov\Delta^2]=C(\ov\Delta^2),$$
with singular set $Z=\ov\Delta^2\cap\{z_1=0\}$. Here $f=(\bar z_1^3,\bar z_2)$ and $\det B_f=3\bar z_1^2$.

Since $\bar z_1^3$ vanishes to order $3>2$ at $z_1=0$, it can be approximated in $C^2$ by functions vanishing near $z_1=0$ \cite[Ex.~6.2]{HW}. Explicitly, take $\chi\in C^\infty(\mathbb R)$ with $\chi=0$ on $(-\infty,1]$ and $\chi=1$ on $[2,\infty)$, let $\delta_l\to0$, and put
$$\phi_l(z_1)=\chi(|z_1|/\delta_l)\,\bar z_1^3,\qquad g_l=(\phi_l(z_1),\bar z_2).$$
The difference $\bar z_1^3-\phi_l=\big(1-\chi(|z_1|/\delta_l)\big)\bar z_1^3$ is supported in $\{|z_1|\le2\delta_l\}$. The factor $1-\chi(|z_1|/\delta_l)$ is bounded, and its derivatives of order $i\ge1$ vanish outside $\{\delta_l\le|z_1|\le2\delta_l\}$. There the derivatives of $|z_1|$ of order $p\ge1$ are $O(\delta_l^{1-p})$, so the derivatives of order $i$ of the factor are $O(\delta_l^{-i})$. The derivatives of $\bar z_1^3$ of order $j-i$ are $O(\delta_l^{3-j+i})$ on $\{|z_1|\le2\delta_l\}$. By the Leibniz rule every derivative of order $j\le2$ of the difference is $O(\delta_l^{3-j})$, so $g_l\to f$ in $C^2$.

On $\{|z_1|<\delta_l/2\}\times\ov\Delta$ every element of $P[z,g_l;\ov\Delta^2]$ is a uniform limit of polynomials in $z_1,z_2,\bar z_2$, hence holomorphic in $z_1$ for fixed $z_2$. So $\bar z_1\notin P[z,g_l;\ov\Delta^2]$, and a fortiori $P[z^{[m]},g_l;\ov\Delta^2]\ne C(\ov\Delta^2)$.
\end{Remark}

\section{Some remarks and  questions}
\begin{Remark}
 Duval's theorem \cite{Du} states: if $f$ is $C^1$ near the closed unit disk $\ov\Delta$ and $|f_{\bar z}|>|f_z|$ on $\ov\Delta$, then the graph of $f$ over $\ov\Delta$ is polynomially convex. Since $f_{\bar z}\ne0$, it follows that $P[z,f;\ov\Delta]=C(\ov\Delta)$. The proof uses Oka's criterion: through each point off the graph one finds an algebraic curve avoiding the graph and deforms it to infinity. The curves are the quadrics $\{(z-a)(w-f(a))=\alpha f_{\bar z}(a)\}$, with $\alpha$ in a small half-disk $\{|\alpha|\le\ve,\ \Re\alpha<0\}$. The key estimate \cite[Lemme~1]{Du} is
$$\fr{(z-a)(f(z)-f(a))}{f_{\bar z}(a)}=\fr{f_z(a)}{f_{\bar z}(a)}(z-a)^2+|z-a|^2\big(1+\varphi(z,a)\big),\qquad\varphi\to0\text{ on the diagonal},$$
which gives a non-negative real part for $|z-a|\le\delta$. This is the positivity required in Proposition \ref{prop:peak}, for $Q_a(z,w)=(z-a)(w-f(a))/f_{\bar z}(a)$, but only for $z$ near $a$. Duval passes from local to global by the injectivity of $f$ \cite[\S1]{Du}, or, in general, by the finiteness of the fibres of $f$, Lagrange interpolation and a deformation of the curves \cite[\S2]{Du}. Proposition \ref{prop:peak} and Theorem \ref{thmMon} are the cases in which the positivity holds globally; then no deformation is needed. The local part carries over to all dimensions.
\end{Remark}
\begin{Proposition}[Duval's estimate in $\C^N$]\label{prop:locDuval}
Let $L\subset\C^N$ be compact and convex, $f\in C^1$ near $L$ with values in $\C^N$, and $U\colon L\to GL_N(\C)$ continuous with
\begin{equation}\label{eq:DuvN}
\Re\big(v^TU(a)B_f(a)\bar v\big)>\big|v^TU(a)A_f(a)v\big|\qquad(a\in L,\ v\in\C^N\setminus\{0\}).
\end{equation}
Then there is $\delta>0$ such that for every compact $K\subset L$ with $\operatorname{diam}K\le\delta$, the graph $\Gam$ is polynomially convex and $P[z,f;K]=C(K)$. If moreover $K$ satisfies (ii), then $P[z^{[m]},f;K]=C(K)$.

If $B_f$ is invertible on $L$, one may take $U=B_f^{-1}$, and \eqref{eq:DuvN} becomes $|v^TB_f^{-1}A_fv|<|v|^2$; this holds if $\|A_f\|<\sigma_{\min}(B_f)$ on $L$. For $N=1$ it is Duval's condition $|f_z|<|f_{\bar z}|$.
\end{Proposition}
\begin{proof}
\emph{Constants.} The function $(a,v)\mapsto\Re(v^TU(a)B_f(a)\bar v)-|v^TU(a)A_f(a)v|$ is continuous and positive on the compact set $L\times\{|v|=1\}$. Let $c>0$ be its minimum and $M=\max_L\|U\|$. Since $A_f,B_f$ are uniformly continuous on $L$, the non-decreasing function
$$\omega(t)=\sup\big\{\|A_f(x)-A_f(y)\|+\|B_f(x)-B_f(y)\|:\ x,y\in L,\ |x-y|\le t\big\}$$
tends to $0$ as $t\to0$. Choose $\delta>0$ with $M\omega(\delta)\le c/2$.

\emph{Taylor estimate.} Let $a,z\in L$, $v=z-a$, $\gamma(t)=a+tv\in L$. Then
$$f(z)-f(a)-A_f(a)v-B_f(a)\bar v=\int_0^1\Big(\big(A_f(\gamma)-A_f(a)\big)v+\big(B_f(\gamma)-B_f(a)\big)\bar v\Big)dt,$$
whose norm is at most $\omega(|v|)|v|$, as $|\gamma(t)-a|\le|v|$. Hence, if $|v|\le\delta$,
\begin{align*}
\Re\Big(v^TU(a)\big(f(z)-f(a)\big)\Big)
&\ge\Re\big(v^TU(a)B_f(a)\bar v\big)-\big|v^TU(a)A_f(a)v\big|-M\omega(|v|)|v|^2\\
&\ge c|v|^2-\tfrac c2|v|^2=\tfrac c2|v|^2 .
\end{align*}

\emph{Polynomial convexity.} Let $K\subset L$ be compact with $\operatorname{diam}K\le\delta$. Then $\widehat K\subset\conv K\subset L$, and $|z-a|\le\operatorname{diam}\conv K=\operatorname{diam}K\le\delta$ for $z\in K$, $a\in\widehat K$. For $a\in\widehat K$ put $Q_a(z,w)=(z-a)^TU(a)(w-f(a))$, a polynomial with $Q_a(a,w)=0$. By the estimate, $\Re Q_a(z,f(z))\ge\frac c2|z-a|^2>0$ for $z\in K\setminus\{a\}$. Proposition \ref{prop:peakhull} (with $W=\C^N$) shows that $\Gam$ is polynomially convex. Note that $K$ itself need not be polynomially convex.

\emph{Density.} If $B_f(a)\bar v=0$ for some $v\ne0$, then \eqref{eq:DuvN} gives $0>|v^TU(a)A_f(a)v|$, which is absurd. So $B_f$ is invertible on $L$, hence on a neighbourhood $V$ of $L$. By Lemma \ref{lem:tr}, the graph of $f$ over $V$ is a totally real $C^1$ manifold containing $\Gam$. Theorem \ref{thm:Iz} with $E=\emptyset$ and \eqref{eq:GamK} give $P[z,f;K]=C(K)$. Theorem \ref{thmA} gives the last claim.

\emph{The case $U=B_f^{-1}$.} This $U$ is continuous on $L$. Then $v^TU(a)B_f(a)\bar v=v^T\bar v=|v|^2$, and
$$|v^TB_f^{-1}A_fv|\le\|B_f^{-1}A_f\|\,|v|^2\le\fr{\|A_f\|}{\sigma_{\min}(B_f)}\,|v|^2 .$$
For $N=1$, \eqref{eq:DuvN} with $U=1/f_{\bar z}$ reads $|f_z/f_{\bar z}|<1$.
\end{proof}

\begin{Remark}[What is known]\label{rem:known}
(1) The polynomial convexity of small compact subsets of a totally real manifold is classical, also for $C^1$ manifolds. Every point of a totally real $C^1$ manifold has a compact neighbourhood in the manifold that is polynomially convex and carries $P=C$ \cite[Cor.~1.6.15, p.~60]{St}; hence its compact subsets are polynomially convex as well. Weinstock's local result for graphs is \cite[Cor.~1.6.17]{St}. The proof in \cite[p.~60]{St} uses a quadratic polynomial of the same kind as $Q_a$. For totally real $C^1$ manifolds there is also a non-negative strictly plurisubharmonic function vanishing exactly on the manifold (Chirka, \cite[Thm.~6.1.6]{St}).

(2) In the $C^2$ case the local statement can be seen directly. Let $\Sigma$ be a $C^2$ submanifold without complex tangents. The square $d^2$ of the distance to $\Sigma$ is $C^2$ and strictly plurisubharmonic near $\Sigma$ \cite[Lemma~17.2]{AW}. Let $B$ be a ball on which $d^2$ is $C^2$ with Levi form $\ge\lambda|v|^2$, $\lambda>0$, let $X\subset\Sigma\cap B$ be compact, and let $p\in\widehat X$. Then $p\in\conv X\subset B$, and Theorem \ref{thm:psh} gives $d^2(p)\le\sup_Xd^2=0$, so $p\in\Sigma$. If $p\notin X$, then $d^2-c|z-p|^2$ is plurisubharmonic on $B$ for $0<c\le\lambda$. It vanishes at $p$ and is negative on $X$, which contradicts Theorem \ref{thm:psh}.

(3) So the conclusion of Proposition \ref{prop:locDuval} is known in substance. What the proposition adds is the mechanism: Duval's estimate in class $C^1$, with an explicit $\delta$ depending only on $c$, $\max\|U\|$ and the modulus of continuity of $df$. It also shows that for $N\ge2$ the missing ingredient is the passage from local to global, which Duval achieves by deforming curves. If $U$ is constant and \eqref{eq:DuvN} holds on all of $L$, the integral along segments gives Proposition \ref{prop:D} (and Remark \ref{rem:shift} with $S=0$) for $K$ of arbitrary diameter. Proposition \ref{prop:prod} is the product case of Duval's theorem.
\end{Remark}

\begin{Problem}
\par (1) Let $f\in C^1$ near the closed ball (or polydisk) $L$ in $\C^N$ satisfy \eqref{eq:DuvN} with $U=B_f^{-1}$, for instance $\|A_f\|<\sigma_{\min}(B_f)$ on $L$. By Proposition \ref{prop:locDuval}, graphs over compact subsets of small diameter are polynomially convex. Is the graph over $L$ polynomially convex? For $N=1$ this is Duval's theorem.
\par (2) Is (iii) necessary in some form when $f$ is real-analytic?
\par (3)
For $f=\bar z+R$ with $R$ Lipschitz with constant $k<1$ on $K$, no condition on the set where $\det B_f=0$ is needed, by Theorem \ref{thmLip}. Can the condition on $Z$ be removed in Corollary \ref{cor:strat} for general $f$ with polynomially convex graph, when $B_f$ is invertible on $\ov{\operatorname{int}K}$? Remark \ref{rem:fibre} shows that some condition is needed when $B_f$ degenerates on the interior.
\par (4)
In Remark \ref{rem:sharp} the failure comes from separation. Assume that $z^{[m]},f$ separate the points of $K$, that $P[z,f;K]=C(K)$, and that \eqref{eq:sector} holds with a larger angle. Can density fail?
\end{Problem}

\section{Declaration of generative AI use}
The research questions, mathematical ideas, and direction of this work were determined
by the authors. Claude's AI was used to assist with literature searches, checks
of mathematical arguments and examples, LaTeX preparation, and language editing. The
authors reviewed the AI-assisted material and take full responsibility for the mathematical
statements, proofs, references, and final manuscript.

\end{document}